\documentclass[a4paper,oneside,12pt]{article}
\usepackage{amsmath}
\usepackage{amssymb}
\usepackage{amsthm}
\usepackage[a4paper,margin=3cm,centering,nohead,ignorefoot,footskip=5ex]{geometry}
\usepackage{url}
\theoremstyle{plain}
\newtheorem{theorem}{Theorem}
\newtheorem{proposition}[theorem]{Proposition}
\newtheorem{lemma}[theorem]{Lemma}
\newtheorem{application}[theorem]{Application}
\theoremstyle{definition}
\newtheorem{definition}[theorem]{Definition}
\newcommand*{\st}{\mathrm{st}}
\newcommand*{\EHAw}{\mathsf{E\text{-}HA}^\omega}
\newcommand*{\EPAw}{\mathsf{E\text{-}PA}^\omega}
\newcommand*{\EHAwst}{\mathsf{E\text{-}HA}_\st^\omega}
\newcommand*{\EPAwst}{\mathsf{E\text{-}PA}_\st^\omega}
\newcommand*{\HA}{\mathsf{HA}}
\newcommand*{\HAwbd}{\mathsf{HA}^\omega_\unlhd}
\newcommand*{\EHAwsst}{\mathsf E\text{-}\mathsf{HA}_\st^{\omega*}}
\newcommand*{\mACw}{\mathsf{mAC}^\omega}
\newcommand*{\Rw}{\mathsf R^\omega}
\newcommand*{\IPwmexistsstf}{\mathsf{IP}_{\mnexistsst}^\omega}
\newcommand*{\MAJw}{\mathsf{MAJ}^\omega}
\newcommand*{\Iw}{\mathsf I^\omega}
\newcommand*{\IPmforallwst}{\mathsf{IP}_{\mforallst}^\omega}
\newcommand*{\Mw}{\mathsf M^\omega}
\newcommand*{\BUDw}{\mathsf{BUD}^\omega}
\newcommand*{\Tforall}{\mathsf T_{\!\forall}}
\newcommand*{\Texists}{\mathsf T_{\!\exists}}
\newcommand*{\USzero}{\mathsf{US_0}}
\newcommand*{\R}{\mathsf R}
\newcommand*{\HGMPst}{\mathsf{HGMP}^\st}
\newcommand*{\PP}{\mathsf P}
\newcommand*{\LEM}{\mathsf{LEM}}
\newcommand*{\Forall}[1]{\forall#1\,}
\newcommand*{\Exists}[1]{\exists#1\,}
\newcommand*{\Finvv}[1]{\Finv#1\,}
\newcommand*{\forallst}{\forall^\mathrm{st}}
\newcommand*{\existsst}{\exists^\mathrm{st}}
\newcommand*{\Forallst}[1]{\forall^\mathrm{st}#1\,}
\newcommand*{\Existsst}[1]{\exists^\mathrm{st}#1\,}
\newcommand*{\Finvvst}[1]{\Finv^\mathrm{st}#1\,}
\newcommand*{\mforall}{\tilde\forall}
\newcommand*{\mexists}{\tilde\exists}
\newcommand*{\mForall}[1]{\tilde\forall#1\,}
\newcommand*{\mExists}[1]{\tilde\exists#1\,}
\newcommand*{\mforallst}{\tilde\forall^\mathrm{st}}
\newcommand*{\mexistsst}{\tilde\exists^\mathrm{st}}
\newcommand*{\mForallst}[1]{\tilde\forall^\mathrm{st}#1\,}
\newcommand*{\mExistsst}[1]{\tilde\exists^\mathrm{st}#1\,}
\newcommand*{\mnexistsst}{{\tilde\nexists^\mathrm{st}}}
\newcommand*{\bb}{\mathrm b}
\newcommand*{\B}{\mathrm B}
\newcommand*{\BB}{{\mathrm B'}}
\newcommand*{\ul}[1]{\mspace{0.75mu}\underline{\mspace{-0.75mu}#1\mspace{-0.75mu}}\mspace{0.75mu}}
\newcommand*{\defeq}{\mathrel{\mathop:}=}
\newcommand*{\defequiv}{\mathrel{\mathop:}\equiv}
\newcommand*{\la}[2]{\lambda#1\,.\,#2}

\begin{document}

\title{Intuitionistic nonstandard bounded modified realisability and functional interpretation\thanks{The authors are grateful to Fernando Ferreira for his availability for discussion and useful remarks and to an anonymous referee for his remarks and suggestions. The authors are however the only persons responsible for any mistakes. 2010 Mathematics Subject Classification: 03B15, 03B20, 03F10, 03F25, 03F35, 03F50, 11U10. Keywords: intuitionism, bounded functional interpretation, bounded modified realisability, majorisability, nonstandard arithmetic, transfer principle.}}
\author{Bruno Dinis\thanks{Departamento de Matem\'atica, Faculdade de Ci\^encias da Universidade de Lisboa, Campo Grande, Edif\'icio~C6, 1749-016~Lisboa, Portugal. \protect\url{bmdinis@fc.ul.pt}. Funded by the grant~SFRH/BPD/97436/2013 of the Funda\c c\~ao para a Ci\^encia e a Tecnologia and by the grant~UID/MAT/04561/2013 of the Centro de Matem\'atica, Aplica\c c\~oes Fundamentais e Investiga\c c\~ao Operacional~/ Funda\c c\~ao da Faculdade de Ci\^encias da Universidade de Lisboa.}
\and
Jaime Gaspar\thanks{School of Computing, University of Kent, Canterbury, Kent, CT2~7NF, United Kingdom. Centro de Matem\'atica e Aplica\c c\~oes (CMA), FCT, UNL. \protect\url{mail@jaimegaspar.com}~/ \protect\url{jg478@kent.ac.uk}, \protect\url{www.jaimegaspar.com}~/ \protect\url{www.cs.kent.ac.uk/people/rpg/jg478}. Funded by a Research Postgraduate Scholarship from the Engineering and Physical Sciences Research Council~/ School of Computing, University of Kent.}}
\date{10 December 2017}
\maketitle

\begin{abstract}
  We present a bounded modified realisability and a bounded functional interpretation of intuitionistic nonstandard arithmetic with nonstandard principles.

  The functional interpretation is the intuitionistic counterpart of Ferreira and Gaspar's functional interpretation and has similarities with Van den Berg, Briseid and Safarik's functional interpretation but replacing finiteness by majorisability.

  We give a threefold contribution: constructive content and proof-theoretical properties of nonstandard arithmetic; filling a gap in the literature; being in line with nonstandard methods to analyse compactness arguments.
\end{abstract}

\section{Introduction}

In the past few years there has been a growing interest in the study of nonstandard arithmetic using realisabilities and functional interpretations. Particularly interesting for us are:
\begin{enumerate}
  \item Ferreira and Gaspar's functional interpretation~\cite{FerreiraGaspar2015} (which deals with majorisability);
  \item Van den Berg, Briseid and Safarik's realisability and functional interpretation~\cite{BergBriseidSafarik2012} (which deals with finiteness).
\end{enumerate}
In their spirit we present
\begin{enumerate}
  \item a realisability (based on Ferreira and Nunes's bounded modified realisability~\cite{FerreiraNunes2006}, and similar to Van den Berg, Briseid and Safarik's realisability~\cite{BergBriseidSafarik2012});
  \item a functional interpretation (based on Ferreira and Oliva's functional interpretation~\cite{FerreiraOliva2005}, and Ferreira and Gaspar's functional interpretation~\cite{FerreiraGaspar2015}, and similar to Van den Berg, Briseid and Safarik's functional interpretation~\cite{BergBriseidSafarik2012});
\end{enumerate}
of intuitionistic nonstandard arithmetic with nonstandard principles and we prove their soundness and characterisation theorems.

Overall we give a threefold contribution:
\begin{enumerate}
  \item we give as applications various results, for example about
  \begin{enumerate}
    \item the constructive content (such as bounded-program extraction and extraction of bounds on witnesses);
     \item proof-theoretical properties (such as relative consistency results and independence results);
  \end{enumerate}
  of nonstandard arithmetic;
  \item we fill an existing gap in the literature, for example
  \begin{enumerate}
    \item we show that if we restrict our realisability to the ``purely external fragment'' (where there are only quantifiers of the form $\existsst$ and $\forallst$), then we recover Ferreira and Nunes's bounded modified realisability, and analogously if we restrict our functional interpretation to the ``purely external fragment'', then we recover Ferreira and Oliva's bounded functional interpretation;
    \item we think that if we combine our intuitionistic functional interpretation with a suitable negative translation, then we obtain Ferreira and Gaspar's classical functional interpretation (a result that we hope to publish soon);
  \end{enumerate}
  \item we are in line with the argument that nonstandard methods can be used to analyse compactness arguments~\cite{FerreiraOliva2005,FerreiraNunes2006,FerreiraGaspar2015,DinisFerreira2017}.
\end{enumerate}

Our realisability and functional interpretation, like previous ones, make use of Nelson's syntactic approach to nonstandard analysis called internal set theory~\cite{Nelson1977, Nelson1988} by extending the language with a new unary predicate $\st(t)$ (meaning ``$t$ is standard'').

\section{Framework}

Let $\EPAw$ and $\EHAw$ be (respectively) Peano and Heyting arithmetics in all finite types with full extensionality and with primitive equality only at type~$0$. In the next definition, proposition and theorem, we recall well-known facts about the Howard-Bezem strong majorisability~$\leq^*_\sigma$~\cite{Bezem1985,Howard1973}.

\begin{definition}\mbox{}
  \begin{enumerate}
    \item The \emph{Howard-Bezem strong majorisability}~$\leq^*_\sigma$ is defined by recursion on the finite type~$\sigma$ by:
    \begin{enumerate}
      \item $s \leq^*_0 t \defequiv s \leq_0 t$;
      \item $s \leq^*_{\rho \rightarrow \sigma} t \defequiv \Forall v \Forall{u \leq^*_\rho v} (su \leq^*_\sigma tv \wedge tu \leq^*_\sigma tv)$.
    \end{enumerate}
    \item We say that $x^\sigma$ is \emph{monotone} if and only if $x \leq^*_\sigma x$.
  \end{enumerate}
\end{definition}

\begin{proposition}
  We have:
  \begin{enumerate}
    \item $\EHAw \vdash x \leq^*_\sigma y \rightarrow y \leq^*_\sigma y$;
    \item $\EHAw \vdash x \leq^*_\sigma y \wedge y \leq^*_\sigma z \rightarrow x \leq^*_\sigma z$.
  \end{enumerate}
\end{proposition}

\begin{theorem}[{Howard's majorisability theorem~\cite{Howard1973}}]
  For all \emph{closed} terms $t^\sigma$ of $\EHAw$, there is a closed term $\tilde t^\sigma$ of $\EHAw$ such that $\EHAw \vdash t \leq^*_\sigma \tilde t$.
\end{theorem}

We introduce a nonstandard variant~$\EHAwst$ of $\EHAw$ analogously to the nonstandard variant~$\EPAwst$~\cite{FerreiraGaspar2015} of $\EPAw$.

\begin{definition}
  The \emph{nonstandard Heyting arithmetic in all finite types with full extensionality}~$\EHAwst$ is obtained from $\EHAw$ by enriching the language and the axioms of $\EHAw$ as follows.
\begin{enumerate}
  \item We add the \emph{standard predicates} $\st^\sigma(t^\sigma)$ (for each finite type $\sigma$).
  \item We add the \emph{standardness axioms}:
    \begin{enumerate}
      \item $x =_\sigma y \wedge \st^\sigma(x) \rightarrow \st^\sigma(y)$;
      \item $\st^\sigma(y) \wedge x \leq_\sigma^* y \rightarrow \st^\sigma(x)$;
      \item $\st^\sigma(t)$ for each closed term $t$;
      \item $\st^{\sigma \rightarrow \tau}(x) \wedge \st^\sigma(y) \rightarrow \st^{\tau}(xy)$;
    \end{enumerate}
    (for each finite types $\sigma$ and $\tau$).
    \item We add the \emph{external induction rule} $\frac{\Phi(0) \ \ \Forall{x^0} (\st^0(x) \to (\Phi(x) \to \Phi(x + 1)))}{\Forall{x^0} (\st^0(x) \to \Phi(x))}$ (or, equivalently, the corresponding axioms).
    \item The logical axioms are extended to the formulas in the language of $\EHAwst$ but the arithmetical axioms are not.
  \end{enumerate}
\end{definition}

The standardness axioms are (essentially) the same as used by Ferreira and Gaspar~\cite{FerreiraGaspar2015}, and by Van den Berg, Briseid and Safarik~\cite{BergBriseidSafarik2012}, with the exception of the second one. The second standardness axiom has a clear meaning in type 0 (namely that the nonstandard natural numbers are an end-extension of the standard natural numbers~\cite[page~12]{Kaye1991}) but not in higher types.

We will use the following convenient abbreviations.
\begin{definition}
  Let $\ul x$ be a tuple of variables $x_1,\ldots,x_n$ of $\EHAwst$ of length $n$, $\ul s$ and $\ul t$ be tuples of terms $s_1,\ldots,s_n$ and $t_1,\ldots,t_n$ (respectively) of $\EHAwst$ also of length $n$, and $\Phi(\ul x)$ and $\Psi(\ul x)$ be formulas of $\EHAwst$. We abbreviate:
  \begin{enumerate}
    \item $s \leq^*_\sigma t$ by $s \leq^*\! t$;
    \item $s_1 \leq^*\! t_1 \wedge \cdots \wedge s_n \leq^*\! t_n$ by $\ul s \leq^*\! \ul t$;
    \item $\st^\sigma(t)$ by $\st(t)$;
    \item $\st(t_1) \wedge \cdots \wedge \st(t_n)$ by $\st(\ul t)$;
    \item $\Forall{\ul x} \Phi(\ul x)$ by $\Forall{\ul x} \Phi$, and analogously for $\Exists{\ul x} \Phi(\ul x)$;
    \item $\Forall{\ul x} (\Phi \to \Psi)$ with
      \begin{gather*}
        \Phi \equiv \ul x \leq^*\! \ul x, \quad
        \Phi \equiv \st(\ul x), \quad
        \Phi \equiv \ul x \leq^*\! \ul t, \quad
        \Phi \equiv \ul x \leq^*\! \ul x \wedge \st(\ul x), \\
        \Phi \equiv \ul x \leq^*\! \ul x \wedge \ul x \leq^*\! \ul t, \quad
        \Phi \equiv \st(\ul x) \wedge \ul x \leq^*\! \ul t, \quad
        \Phi \equiv \ul x \leq^*\! \ul x \wedge \st(\ul x) \wedge \ul x \leq^*\! \ul t
      \end{gather*}
    by (respectively)
    \begin{equation*}
      \mForall{\ul x} \Psi, \ \ 
      \Forallst{\ul x} \Psi, \ \ 
      \Forall{\ul x \leq^*\! \ul t} \Psi, \ \ 
      \mForallst{\ul x} \Psi, \ \ 
      \mForall{\ul x \leq^*\! \ul t} \Psi, \ \ 
      \Forallst{\ul x \leq^*\! \ul t} \Psi, \ \ 
      \mForallst{\ul x \leq^*\! \ul t} \Psi,
    \end{equation*}
    and analogously for $\Exists{\ul x} (\Phi \wedge \Psi)$.
   \end{enumerate}
\end{definition}

Many times we have a formula $\Forall{\ul y} \Exists{\ul z} \Phi(\ul z)$ and we take $\ul z$ as functions of $\ul y$ obtaining a new formula; then we denote those functions by $\ul Z$ (an uppercase italic letter), so the new formula is denoted by $\Exists{\ul Z} \Forall{\ul y} \Phi(\ul Z\ul y)$. Sometimes we have a formula $\Forall{\ul x} \Exists{\ul Z} \Forall{\ul y} \Phi(\ul Z\ul y)$ and we take $\ul Z$ as functions of $\ul x$ obtaining a new formula; then we denote those functions by $\ul{\mathrm Z}$ (an uppercase roman letter), so the new formula is denoted by $\Exists{\ul{\mathrm Z}} \Forall{\ul x,\ul y} \Phi(\ul{\mathrm Z}\ul x\ul y)$. We introduce this convention in the next definition.

\begin{definition}
  Let $\ul x$, $\ul y$ and $\ul z$ be tuples of variables.
  \begin{enumerate}
    \item If the variables $\ul z$ become functions of variables $\ul y$, then we denote those functions by $\ul Z$;
    \item If the variables $\ul Z$ become functions of variables $\ul x$, then we denote those functions by $\ul{\mathrm Z}$.
  \end{enumerate}
\end{definition}

We follow the following conventions by Nelson~\cite{Nelson1977}.
\begin{definition}\mbox{}
  \begin{enumerate}
    \item We say that a formula of $\EHAwst$ is \emph{internal} if and only if the standard predicates do not occur in the formula and we denote the formula by lower case Greek letters $\phi,\psi,\ldots$\:.
    \item We say that a formula of $\EHAwst$ is \emph{external} if and only if the standard predicates may occur in the formula and we denote the formula by upper case Greek letters $\Phi,\Psi,\ldots$\:.
  \end{enumerate}
\end{definition}

\section{Intuitionistic nonstandard bounded modified realisability}

\begin{definition}
  The \emph{intuitionistic nonstandard bounded modified realisability}~$\bb$ assigns to each formula $\Phi$ of $\EHAwst$ the formulas $\Phi^\bb$ and $\Phi_\bb(\ul a)$ of $\EHAwst$ such that $\Phi^\bb \equiv \mExistsst{\ul a} \Phi_\bb(\ul a)$ accordingly to the following clauses where $\Phi_\bb(\ul a)$ is the part inside square brackets below. For atomic formulas, we define:
  \begin{enumerate}
    \item $\Phi^\bb \defequiv [\Phi]$ for internal atomic formulas $\Phi$ (so the tuple $\ul a$ is empty);
    \item $\st(t)^\bb \defequiv\mExistsst a [t \leq^*\! a]$;
  \end{enumerate}
  For the remaining formulas, if $\Phi^\bb \equiv \mExistsst {\ul a} \Phi_\bb(\ul a)$ and $\Psi^\bb \equiv \mExistsst{\ul b} \Psi_\bb(\ul b)$, then we define:
  \begin{enumerate}
    \setcounter{enumi}{2}
    \item $(\Phi \wedge \Psi)^\bb \defequiv \mExistsst{\ul a,\ul b} [\Phi_\bb(\ul a) \wedge \Psi_\bb(\ul b)]$;
    \item $(\Phi \vee \Psi)^\bb \defequiv \mExistsst{\ul a,\ul b} [\Phi_\bb(\ul a) \vee \Psi_\bb(\ul b)]$;
    \item $(\Phi \rightarrow \Psi)^\bb \defequiv \mExistsst{\ul B} [\mForallst{\ul a} (\Phi_\bb(\ul a) \rightarrow \Psi_\bb(\ul B\ul a))]$;
    \item $(\Forall x \Phi)^\bb \defequiv \mExistsst{\ul a} [\Forall x \Phi_\bb(\ul a)]$;
    \item $(\Exists x \Phi)^\bb \defequiv \mExistsst{\ul a} [\Exists x \Phi_\bb(\ul a)]$.
  \end{enumerate}
\end{definition}

We need the next lemmas to prove the soundness and characterisation theorems of $\bb$ and to do the applications of $\bb$.

\begin{lemma}\label{lemma realisability of defined quantifications}\mbox{}
  \begin{enumerate}
    \item For all internal formulas~$\phi$ of $\EHAwst$, we have $\phi^\bb \equiv \phi_\bb(\ul a) \equiv \phi$ (so the tuple $\ul a$ is empty).
    \item For all formulas~$\Phi$ of $\EHAwst$, we have (the equivalences being provable in $\EHAwst$):
    \begin{enumerate}
     \item $(\Forallst{\ul x} \Phi)^\bb \leftrightarrow \mExistsst{\ul A} [\mForallst{\ul b} \Forall{\ul x \leq^*\! \ul b} \Phi_\bb(\ul A\ul b)]$;
      \item\label{item realisability of defined quantifications} $(\mForallst{\ul x} \Phi)^\bb \leftrightarrow \mExistsst{\ul A} [\mForallst{\ul b} \mForall{\ul x \leq^*\! \ul b} \Phi_\bb(\ul A\ul b)]$;
      \item $(\Existsst{\ul x} \Phi)^\bb \equiv \mExistsst{b,\ul a} [\Exists{\ul x \leq^*\! \ul b} \Phi_\bb(\ul a)]$;
      \item $(\mExistsst{\ul x} \Phi)^\bb \equiv \mExistsst{b,\ul a} [\mExists{\ul x \leq^*\! \ul b} \Phi_\bb(\ul a)]$;
      \item $(\Exists{\ul x \leq^*\! \ul t} \Phi)^\bb \equiv \mExistsst{\ul a} [\Exists{\ul x \leq^*\! \ul t} \Phi_\bb(\ul a)]$;
      \item $(\mExists{\ul x \leq^*\! \ul t} \Phi)^\bb \equiv \mExistsst{\ul a} [\mExists{\ul x \leq^*\! \ul t} \Phi_\bb(\ul a)]$.
    \end{enumerate}
  \end{enumerate}
\end{lemma}

\begin{proof}\mbox{}
  \begin{enumerate}
    \item The proof is by a simple induction on the length of $\phi$ (analogous to the proof of lemma~\ref{lemma realisability of mexistsst-free formula} below).
    \item The proof is by simple calculations using that $\ul x \leq^*\! \ul t$ is an internal formula. For example, for point~\ref{item realisability of defined quantifications}:
    \begin{align*}
      \Phi^\bb &\equiv \mExistsst{\ul a} \Phi_\bb(\ul a), &\text{(assumption)} \\
      (\ul x \leq^*\! \ul x)^\bb &\equiv \ul x \leq^*\! \ul x, &\text{($\ul x \leq^*\! \ul x$ internal)} \\
      \st(\ul x)^\bb &\equiv \mExistsst{\ul b} [\ul x \leq^*\! \ul b], &\text{(calculation)} \\
      (\ul x \leq^*\! \ul x \wedge \st(\ul x))^\bb &\equiv \mExistsst{\ul b} [\ul x \leq^*\! \ul x \wedge \ul x \leq^*\! \ul b], &\text{(calculation)} \\
      (\ul x \leq^*\! \ul x \wedge \st(\ul x) \to \Phi)^\bb &\equiv \mExistsst{\ul A} [\mforallst{\ul b} \\ &\phantom{{}\equiv{}} (\ul x \leq^*\! \ul x \wedge \ul x \leq^*\! \ul b \to \Phi_\bb(\ul A\ul b))], &\text{(calculation)} \\
      (\Forall{\ul x} (\ul x \leq^*\! \ul x \wedge \st(\ul x) \to \Phi))^\bb &\equiv \mExistsst{\ul A} [\Forall{\ul x} \mforallst{\ul b} \\ &\phantom{{}\equiv{}} (\ul x \leq^*\! \ul x \wedge \ul x \leq^*\! \ul b \to \Phi_\bb(\ul A\ul b))], &\text{(calculation)} \\
      (\Forall{\ul x} (\ul x \leq^*\! \ul x \wedge \st(\ul x) \to \Phi))^\bb &\leftrightarrow \mExistsst{\ul A} [\mForallst{\ul b} \mForall{\ul x \leq^*\! \ul b} \Phi_\bb(\ul A\ul b)], &\text{(logic)}
    \end{align*}
    where $\mForallst{\ul x} \Phi \equiv \Forall{\ul x} (\ul x \leq^*\! \ul x \wedge \st(\ul x) \to \Phi)$.\qedhere
  \end{enumerate}
\end{proof}

The formulas $\Phi_\bb(\ul a)$ are monotone on $\ul a$ in the sense of the following lemma.

\begin{lemma}[monotonicity of $\bb$]
  For all formulas $\Phi$ of $\EHAwst$, we have $\EHAwst \vdash \Phi_\bb(\ul a) \wedge \ul a \leq^*\! \ul{\tilde a} \rightarrow \Phi_\bb(\ul{\tilde a})$.
\end{lemma}

\begin{proof}
  The proof is by a simple induction on the length of $\Phi$. For example, in the case of $\rightarrow$, we have $(\Phi \rightarrow \Psi)_\bb(\ul B) \equiv \mForallst{\ul a} (\Phi_\bb(\ul a) \rightarrow \Psi_\bb(\ul B\ul a))$, we also have $\EHAwst \vdash \Psi_\bb(\ul b) \wedge \ul b \leq^*\! \ul{\tilde b} \rightarrow \Psi_\bb(\ul{\tilde b})$ by induction hypothesis, and we prove $\mForallst{\ul a} (\Phi_\bb(\ul a) \rightarrow \Psi_\bb(\ul B\ul a)) \wedge \ul B \leq^*\! \ul{\tilde B} \rightarrow \mForallst{\ul a} (\Phi_\bb(\ul a) \rightarrow \Psi_\bb(\ul{\tilde B}\ul a))$ by noticing that $\ul B \leq^*\! \ul{\tilde B}$ implies $\ul B \ul a \leq^*\! \ul{\tilde B}\ul a$ (for monotone $\ul a$) and so the induction hypothesis implies $\Psi_\bb(\ul B\ul a) \rightarrow \Psi_\bb(\ul{\tilde B}\ul a)$ (for monotone $\ul a$).
\end{proof}

Now we define the class of $\mexistsst$-free formulas that has a role similar to the one of $\mexists$-free formulas~\cite{FerreiraNunes2006}, which in turn reminds the well-known $\exists$-free formulas (with the difference that $\mexistsst$-free and $\mexists$-free formulas allow disjunctions).

\begin{definition}
  We say that a formula of $\EHAwst$ is \emph{$\mexistsst$-free} if and only if it is built:
  \begin{enumerate}
    \item from atomic internal formulas~$s =_0 t$;
    \item by conjunctions~$\wedge$;
    \item by disjunctions~$\vee$;
    \item by implications~$\to$;
    \item by quantifications~$\forall$ and $\exists$ (so also $\mforall$ and $\mexists$);
    \item by monotone standard universal quantifications~$\mforallst$ (but not $\mexistsst$).
  \end{enumerate}
  and we denote the formula by upper case Greek letters with $\mnexistsst$ as a subscript as in $\Phi_\mnexistsst,\Psi_\mnexistsst,\ldots$\,.
\end{definition}

\begin{lemma}
  \label{lemma realisability of mexistsst-free formula}
  For all $\mexistsst$-free formulas~$\Phi_\mnexistsst$ of $\EHAwst$, we have $(\Phi_\mnexistsst)^\bb \equiv (\Phi_\mnexistsst)_\bb(\ul a)$ (so the tuple~$\ul a$ is empty) and $\EHAwst \vdash (\Phi_\mnexistsst)_\bb \leftrightarrow \Phi_\mnexistsst$.
\end{lemma}

\begin{proof}
  The proof is by a simple induction on the length of $\Phi_\mnexistsst$. For example, in the case of $\mforallst$, we have $(\mForallst x \Phi_\mnexistsst)^\bb \equiv \mExistsst{\ul A} [\mForallst b \mForall{x \leq^*\! b} (\Phi_\mnexistsst)_\bb(\ul Ab)]$ (where by induction hypothesis $\ul A$ is empty), which is $\mForallst b \mForall{x \leq^*\! b} (\Phi_\mnexistsst)_\bb$ (where $(\Phi_\mnexistsst)_\bb \leftrightarrow \Phi_\mnexistsst$ by induction hypothesis), which is equivalent to $\mForallst x \Phi_\mnexistsst$.
\end{proof}

\begin{lemma}
  \label{lemma realisability is mexistsst-free formula}
  For all formulas~$\Phi$ of $\EHAwst$, the formula $\Phi_\bb(\ul a)$ is $\mexistsst$-free.
\end{lemma}

\begin{proof}
  The proof is by a simple induction on the length of $\Phi$. For example, in the case of $\st(t)$, we have $\st(t)_\bb(a) \equiv t \leq^*\! a$, which is $\mexistsst$-free because it is built from formulas of the form $r \leq_0 s$ (which is an abbreviation for an appropriate internal atomic formula) by means of $\forall$, $\to$ and $\wedge$.
\end{proof}

The following four principles, inspired by or even taken from Ferreira, Nunes and Gaspar's articles~\cite{FerreiraNunes2006,FerreiraGaspar2015}, play an important role in the sequel as they are the characteristic principles of $\bb$.

\begin{definition}
  We define the following principles:
  \begin{enumerate}
    \item the \emph{monotone choice}~$\mACw$ is all instances of $\mForallst x \mExistsst y \Phi \rightarrow \mExistsst Y \mForallst x \mExists{y \leq^*\! Yx} \Phi$;
    \item the \emph{realisation}~$\Rw$ is all instances of $\Forall x \Existsst y \Phi \rightarrow \mExistsst z\Forall x \Exists{y \leq^*\! z} \Phi$;
    \item the \emph{independence of premises}~$\IPwmexistsstf$ is all instance of $(\Phi_\mnexistsst \rightarrow \mExistsst x \Psi) \rightarrow \mExistsst y (\Phi_\mnexistsst \rightarrow \mExists{x \leq^*\! y} \Psi)$;
    \item the \emph{majorisability axioms} $\MAJw$ are all instances of $\Forallst x \Existsst y (x \leq^*\! y)$.
  \end{enumerate}
\end{definition}

The principles above generalise to tuples of variables, that is each principle provably implies in $\EHAwst$ its variant where the single variables~$x$, $y$, $Y$ and $z$ are replaced by (respectively) the tuples of variables~$\ul x$, $\ul y$, $\ul Y$ and $\ul z$ (proof sketch: first we generalise from $x$ to $\ul x$ by induction on the length of $\ul x$ and then we generalise from $y$, $Y$ and $z$ to $\ul y$, $\ul Y$ and $\ul z$ by induction on the length of $\ul y$, $\ul Y$ and $\ul z$).

\label{page Rw implies MAJw}The principle~$\Rw$ implies the principle $\MAJw$, that is $\EHAwst + \Rw$ proves all instances of $\MAJw$ (proof sketch: we take any standard $x'$, we take $\Phi$ to be the formula $y = x'$ in $\Rw$, then the premise of $\Rw$ is provable and the conclusion of $\Rw$ implies $\Existsst z (x' \leq^*\! z)$). However, we keep explicitly mentioning $\MAJw$ (even in the presence of $\Rw$) because of its importance.

The following two principles, the first one inspired by Ferreira and Oliva's article~\cite{FerreiraOliva2005}, are used later on in applications of $\bb$.

\begin{definition}
  We define the following principles:
  \begin{enumerate}
    \item \emph{Markov's principle}~$\Mw$ is all instances of $(\mForallst x \phi \rightarrow \psi) \rightarrow \mExistsst y (\mForall{x \leq^*\! y} \phi \rightarrow \psi)$;
    \item the \emph{law of excluded middle}~$\LEM$ is all instances of $\Phi \vee \neg\Phi$.
  \end{enumerate}
\end{definition}

Now we present our main theorem about~$\bb$. This theorem shows that $\bb$:
\begin{enumerate}
  \item interprets $\EHAwst + \mACw + \Rw + \IPwmexistsstf + \MAJw$ into $\EHAwst$;
  \item extracts computational information~$\ul t$ from proofs in $\EHAwst + \mACw + \Rw + \IPwmexistsstf + \MAJw$.
\end{enumerate}

\begin{theorem}[soundness theorem of $\bb$]
  For all formulas $\Phi$ of $\EHAwst$, if
  \begin{equation*}
    \EHAwst + \mACw + \Rw + \IPwmexistsstf + \MAJw \vdash \Phi,
  \end{equation*}
  then there are \emph{closed} monotone terms $\ul t$ of appropriate types such that
  \begin{equation*}
    \EHAwst \vdash \Phi_\bb(\ul t).
  \end{equation*}
\end{theorem}

\begin{proof}
  The proof is by induction on the length of the derivation of $\Phi$. The logical and arithmetical axioms and rules are dealt with in a way similar to the bounded modified realisability~\cite{FerreiraNunes2006}. So let us focus on the standardness axioms and the characteristic principles (except $\MAJw$, which follows from $\Rw$). We will sometimes use implicitly lemmas~\ref{lemma realisability of defined quantifications} and \ref{lemma realisability of mexistsst-free formula}.
  \begin{description}
    \item[$x = y \wedge \st(x) \rightarrow \st(y)$.] Its interpretation is equivalent to $\mExistsst B \mForallst a (x = y \wedge x \leq^*\! a \to y \leq^*\! Ba)$. It is interpreted by $B \defeq \la{a}{a}$.

    \item[$\st(y) \wedge x \leq^*\! y \rightarrow \st(x)$.] Its interpretation is equivalent to $\mExistsst B \mForallst a (y \leq^*\! a \wedge x \leq^*\! y \rightarrow x \leq^*\! Ba)$. It is interpreted by $B \defeq \la{a}{a}$.

    \item[$\st(t)$.] Its interpretation is $\mExistsst a (t \leq^*\! a)$. It is interpreted by a closed term $\tilde t$ such that $t \leq^*\! \tilde t$, which exists by Howard's majorisability theorem.

    \item[$\st(x) \wedge \st(y) \rightarrow \st(xy)$.] Its interpretation is $\mExistsst C\mForallst{a,b} (x \leq^*\! a \wedge y \leq^*\! b \rightarrow xy \leq^*\! Cab)$. It is interpreted by $C \defeq \la{a,b}{ab}$.

    \item[$\mForallst x \mExistsst y \Phi \rightarrow \mExistsst Y \mForallst x \mExists{y \leq^*\! Yx} \Phi$.] Its interpretation is equivalent to
    \begin{equation*}
      \begin{gathered}
        \mExistsst{F,\mathrm D} \mForallst{B,A} (\mForallst c \mForall{x \leq^*\! c} \mExistsst{y \leq^*\! Bc} \Phi_\bb(Ac) \rightarrow \\
        \mExistsst{Y \leq^*\! FBA} \mForallst e \mForall{x \leq^*\! e} \mExists{y \leq^*\! Yx} \Phi_\bb(\mathrm DBAe)).
      \end{gathered}
    \end{equation*}
    It is interpreted by $F \defeq \la{B,A}{B}$ and $\mathrm D \defeq \la{B,A}{A}$. Indeed, from the premise~$\mForallst c \mForall{x \leq^*\! c} \mExistsst{y \leq^*\! Bc} \Phi_\bb(Ac)$ we get $\mForallst x \mExistsst{y \leq^*\! Bx} \Phi_\bb(Ax)$, so (using monotonicity) we get the conclusion~$\mForallst e \mForall{x \leq^*\! e} \mExists{y \leq^*\! Yx} \Phi_\bb(\mathrm DBAe)$ with $Y \defeq B$.

    \item[$\Forall x \Existsst y \Phi \rightarrow \mExistsst z\Forall x \Exists{y \leq^*\! z} \Phi$.] Its interpretation is
    \begin{equation*}
      \mExistsst{D,\ul C} \mForallst{b,\ul a} (\Forall x \Exists{y \leq^*\! b} \Phi_\bb(\ul a) \rightarrow \mExists{z\leq^*\! Db\ul a} \Forall x \Exists{y \leq^*\! z} \Phi_\bb(\ul Cb\ul a)).
    \end{equation*}
    It is interpreted by $D \defeq \la{b,\ul a}{b}$ and $\ul C \defeq \la{b,\ul a}{\ul a}$.

    \item[$(\Phi \rightarrow \mExistsst x \Psi) \rightarrow \mExistsst y (\Phi \rightarrow \mExists{x \leq^*\! y} \Psi)$.] Its interpretation is equivalent to
    \begin{equation*}
      \mExistsst{D,\ul C} \mForallst{b,\ul a} \big((\Phi \rightarrow \mExists{x \leq^*\! b} \Psi_\bb(\ul a)) \rightarrow \mExists{y \leq^*\! Db\ul a} (\Phi \rightarrow \mExists{x \leq^*\!  y} \Psi_\bb(\ul Cb\ul a))\big).
    \end{equation*}
    It is interpreted by $D \defeq \la{b,\ul a}{b}$ and $\ul C \defeq \la{b,\ul a}{\ul a}$.\qedhere
  \end{description}
\end{proof}

\begin{theorem}[characterisation theorem of $\bb$] For all formulas $\Phi$ of $\EHAwst$, we have
  \begin{equation*}
    \EHAwst + \mACw + \Rw + \IPwmexistsstf + \MAJw \vdash \Phi \leftrightarrow \Phi^\bb.
  \end{equation*}
\end{theorem}

\begin{proof}
  The proof is by induction on the length of $\Phi$. The cases of internal atomic formulas, $\wedge$, $\vee$ and $\exists$ are easy, so let us focus on the remaining cases. We sometimes use implicitly lemma~\ref{lemma realisability is mexistsst-free formula}.
  \begin{description}
    \item[$\st(t)$.] We have
    \begin{align*}
      \st(t) &\leftrightarrow &\text{($\MAJw$, $\st(y) \wedge x \leq^*\! y \rightarrow \st(x)$)} \\
      \mExistsst a (t \leq^*\! a) &\equiv &\text{(definition)} \\
      \st(t)^\bb &.
    \end{align*}

    \item[$\rightarrow$.] We have
    \begin{align*}
      (\Phi \rightarrow \Psi) &\leftrightarrow &\text{(induction hypothesis)} \\
      (\mExistsst{\ul a} \Phi_\bb(\ul a) \rightarrow \mExistsst{\ul b} \Psi_\bb(\ul b)) &\leftrightarrow &\text{(logic)} \\
      \mForallst{\ul a} (\Phi_\bb(\ul a) \rightarrow \mExistsst{\ul b} \Psi_\bb(\ul b)) &\leftrightarrow &\text{($\IPwmexistsstf$)} \\
      \mForallst{\ul a} \mExistsst{\ul b} (\Phi_\bb(\ul a) \rightarrow \mExists{\ul{\tilde b} \leq^*\! \ul b} \Psi_\bb(\ul{\tilde b})) &\leftrightarrow &\text{(monotonicity)} \\
      \mForallst{\ul a} \mExistsst{\ul b} (\Phi_\bb(\ul a) \rightarrow \Psi_\bb(\ul b)) &\leftrightarrow &\text{($\mACw$)} \\
      \mExists{\ul B} \mForallst{\ul a} \mExistsst{\ul b \leq^*\! \ul B\ul a} (\Phi_\bb(\ul a) \rightarrow \Psi_\bb(\ul b)) &\leftrightarrow &\text{(monotonicity)} \\
      \mExists{\ul B} \mForallst{\ul a} (\Phi_\bb(\ul a) \rightarrow \Psi_\bb(\ul B\ul a)) &\equiv &\text{(definition)} \\
      (\Phi \rightarrow \Psi)^\bb &.
    \end{align*}

    \item[$\forall$.] We have
    \begin{align*}
      \Forall x \Phi &\leftrightarrow &\text{(induction hypothesis)} \\
      \Forall x \mExistsst{\ul a} \Phi_\bb(\ul a) &\leftrightarrow &\text{($\Rw$)} \\
      \mExistsst{\ul a} \Forall x \mExists{\ul{\tilde a} \leq^*\! \ul a} \Phi_\bb(\ul{\tilde a}) &\leftrightarrow &\text{(monotonicity)} \\
      \mExistsst{\ul a} \Forall x \Phi_\bb(\ul a) &\equiv &\text{(definition)} \\
      (\Forall x \Phi)^\bb &. &&\qedhere
    \end{align*}
  \end{description}
\end{proof}

In the next result we give some applications of $\bb$:
\begin{enumerate}
  \item the first point is an equiconsistency result (``equiconsistency'' meaning that two theories are both consistent or both inconsistent);
  \item the second point is a bounded variant of the existence property (``bounded'' meaning that instead of giving a term~$t$ such that $\Phi(t)$ said to be witnessing an existential quantification~$\Exists x \Phi(x)$, it gives a term~$t$ such that $\Exists{x \leq^*\! t} \Phi(x)$ said to be bounding the existential quantification~$\Exists x \Phi(x)$);
  \item the third, fourth and fifth points are term-and-rule variants of $\mACw$, $\Rw$ and $\IPwmexistsstf$ (``term'' meaning that some existential quantification~$\Exists x \Phi(x)$ is witnessed by some term~$t$, and ``rule'' meaning that instead of $\Phi \to \Psi$ we have $\cdots \vdash \Phi \ \Rightarrow \ \cdots \vdash \Psi$ also denoted by $\frac{\Phi}{\Psi}$);
  \item the sixth point is a conservation result (``conservation'' meaning that if a stronger theory proves a formula of a certain form, then a weaker theory already proves it).
  \item the seventh and eight points are independence results (``independence'' meaning that a formula is neither provable nor refutable by a theory).
\end{enumerate}
Most of these results are inspired by similar folklore results.

\begin{application}\label{realisability applications}
  Let $\EHAwst + \PP \defeq \EHAwst + \mACw + \Rw + \IPwmexistsstf + \MAJw$.
  \begin{enumerate}
    \item The theories $\EHAwst + \PP$ and $\EHAwst$ are equiconsistent.
    \item\label{item realisability application} If $\EHAwst + \PP \vdash \Existsst{\ul x} \Phi$, then there are closed monotone terms~$\ul t$ of appropriate types such that $\EHAwst + \PP \vdash \Existsst{\ul x \leq^*\! \ul t} \Phi$. Moreover, if $\Phi$ is $\mexistsst$-free, then in the conclusion we can replace $\EHAwst + \PP$ by $\EHAwst$. Analogously, if all quantifications of $\ul x$ are monotone.
    \item If $\EHAwst + \PP \vdash \Forallst{\ul x} \Existsst{\ul y} \Phi$, then there are closed monotone terms~$\ul t$ of appropriate types such that $\EHAwst + \PP \vdash \mForallst{\ul z} \Forall{\ul x \leq^*\! \ul z} \Existsst{\ul y \leq^*\! \ul t\ul z} \Phi$. Also, in the conclusion we can replace $\mForallst{\ul z} \Forall{\ul x \leq^*\! \ul z} \Existsst{\ul y \leq^*\! \ul t\ul z} \Phi$ by $\mForallst{\ul x} \Existsst{\ul y \leq^*\! \ul t\ul x} \Phi$. Moreover, if $\Phi$ is $\mexistsst$-free, then in the conclusion we can replace $\EHAwst + \PP$ by $\EHAwst$. Analogously, if all quantifications of $\ul x$ are monotone, or if all quantifications of $\ul y$ are monotone, or if both.
    \item\label{realisability application mACw} If $\EHAwst + \PP \vdash \Forall{\ul x} \Existsst{\ul y} \Phi$, then there are closed monotone terms~$\ul t$ of appropriate types such that $\EHAwst + \PP \vdash \Forall{\ul x} \Existsst{\ul y \leq^*\! \ul t} \Phi$. Moreover, if $\Phi$ is $\mexistsst$-free, then in the conclusion we can replace $\EHAwst + \PP$ by $\EHAwst$. Analogously, if all quantifications of $\ul x$ are monotone, or if all quantifications of $\ul y$ are monotone, or if both.
    \item If $\EHAwst + \PP \vdash \Phi_\mnexistsst \to \Existsst{\ul x} \Psi$, then there are closed monotone terms~$\ul t$ of appropriate types such that $\EHAwst + \PP \vdash \Phi_\mnexistsst \to \Existsst{\ul x \leq^*\! \ul t} \Psi$. Moreover, if $\Psi$ is $\mexistsst$-free, then in the conclusion we can replace $\EHAwst + \PP$ by $\EHAwst$. Analogously, if all quantifications of $\ul x$ are monotone.
    \item If $\EHAwst + \PP \vdash \mForallst{\ul x} \Existsst{\ul y} \Phi_\mnexistsst$, then $\EHAwst \vdash  \mForallst{\ul x} \Existsst{\ul y} \Phi_\mnexistsst$. Analogously, if all quantifications of $\ul y$ are monotone.
    \item If $\EHAwst + \PP$ is consistent, then there is an instance~$\Phi$ of $\Mw$ such that $\EHAwst + \PP \nvdash \Phi$ and $\EHAwst + \PP \nvdash \neg\Phi$.
    \item If $\EHAwst + \PP$ is consistent, then there is an instance~$\Phi$ of $\LEM$ such that $\EHAwst + \PP \nvdash \Phi$ and $\EHAwst + \PP \nvdash \neg \Phi$.
  \end{enumerate}
\end{application}

\begin{proof}
  The ``Moreover, \ldots'' parts are proved analogously to their preceding parts but using lemma~\ref{lemma realisability of mexistsst-free formula} instead of the characterisation theorem to remain in $\EHAwst$ instead of going to $\EHAwst + \PP$. The ``Analogously, \ldots'' parts are proved analogously to their proceedings parts and sometimes are also a corollary to their preceding parts (for example, a proof sketch for point~\ref{item realisability application}: rewrite $\EHAwst + \PP \vdash \mExistsst{\ul x} \tilde\Phi$ as $\EHAwst + \PP \vdash \Existsst{\ul x} (\ul x \leq^*\! \ul x \wedge \tilde\Phi)$, apply the first sentence of point~\ref{item realisability application} with $\Phi \defequiv \ul x \leq^*\! \ul x \wedge \tilde\Phi$ getting $\EHAwst + \PP \vdash \Existsst{\ul x \leq^*\! \ul t} (\ul x \leq^*\! \ul x \wedge \tilde\Phi)$, and then rewrite this as $\EHAwst + \PP \vdash \mExistsst{\ul x \leq^*\! \ul t} \tilde\Phi$). We will sometimes use implicitly that closed terms are standard.
  \begin{enumerate}
    \item We have
    \begin{align*}
      \EHAwst + \PP \vdash 0 =_0 1 &\ \Rightarrow & \text{(soundness)} \\
      \EHAwst \vdash (0 =_0 1)_\bb &\ \Rightarrow & \text{(definition)} \\
      \EHAwst \vdash 0 =_0 1. &
    \end{align*}
    Trivially, $\EHAwst \vdash 0 =_0 1 \ \Rightarrow \ \EHAwst + \PP \vdash 0 =_0 1$.

    \item We have
    \begin{align*}
      \EHAwst + \PP \vdash \Existsst{\ul x} \Phi &\ \Rightarrow &\text{(soundness)} \\
      \EHAwst \vdash (\Exists{\ul x} \Phi)_\bb(\ul t,\ul s) &\ \Rightarrow &\text{(lemma~\ref{lemma realisability of defined quantifications})} \\
      \EHAwst \vdash \Exists{\ul x \leq^*\! \ul t} \Phi_\bb(\ul s) &\ \Rightarrow &\text{($\ul t$ closed)} \\
      \EHAwst \vdash \Existsst{\ul x \leq^*\! \ul t} \Phi_\bb(\ul s) &\ \Rightarrow &\text{($\ul s$ closed monotone)} \\
      \EHAwst \vdash \Existsst{\ul x \leq^*\! \ul t} \Phi^\bb &\ \Rightarrow &\text{(characterisation)} \\
      \EHAwst + \PP \vdash \Existsst{\ul x \leq^*\! \ul t} \Phi &.
    \end{align*}

    \item We have
    \begin{align*}
      \EHAwst + \PP \vdash \Forallst{\ul x} \Existsst{\ul y} \Phi &\ \Rightarrow &\text{(soundness)} \\
      \EHAwst \vdash (\Forallst{\ul x} \Existsst{\ul y} \Phi)_\bb(\ul t,\ul s) &\ \Rightarrow &\text{(lemma~\ref{lemma realisability of defined quantifications})} \\
      \EHAwst \vdash \mForallst{\ul z} \Forall{\ul x \leq^*\! \ul z} \Exists{\ul y \leq^*\! \ul t\ul z} \Phi_\bb(\ul s\ul z) &\ \Rightarrow &\text{($\ul t\ul z$ standard)} \\
      \EHAwst \vdash \mForallst{\ul z} \Forall{\ul x \leq^*\! \ul z} \Existsst{\ul y \leq^*\! \ul t\ul z} \Phi_\bb(\ul s\ul z) &\ \Rightarrow &\text{($\ul s\ul z$ monotone standard)} \\
      \EHAwst \vdash \mForallst{\ul z} \Forall{\ul x \leq^*\! \ul z} \Existsst{\ul y \leq^*\! \ul t\ul z} \Phi^\bb &\ \Rightarrow &\text{(characterisation)} \\
      \EHAwst + \PP \vdash \mForallst{\ul z} \Forall{\ul x \leq^*\! \ul z} \Existsst{\ul y \leq^*\! \ul t\ul z} \Phi &.
    \end{align*}
    Also, $\EHAwst + \PP \vdash \mForallst{\ul z} \Forall{\ul x \leq^*\! \ul z} \Existsst{\ul y \leq^*\! \ul t\ul z} \Phi \ \Rightarrow \ \EHAwst + \PP \vdash \mForallst{\ul x} \Existsst{\ul y \leq^*\! \ul t\ul x} \Phi$ (by taking $\ul x \defeq \ul z$).

    \item We have
    \begin{align*}
      \EHAwst + \PP \vdash \Forall{\ul x} \Existsst{\ul y} \Phi &\ \Rightarrow &\text{(soundness)} \\
      \EHAwst \vdash (\Forall{\ul x} \Existsst{\ul y} \Phi)_\bb(\ul t,\ul s) &\ \Rightarrow &\text{(lemma~\ref{lemma realisability of defined quantifications})} \\
      \EHAwst \vdash \Forall{\ul x} \Exists{\ul y \leq^*\! \ul t} \Phi_\bb(\ul s) &\ \Rightarrow &\text{($\ul t$ closed monotone)} \\
      \EHAwst \vdash \Forall{\ul x} \Existsst{\ul y \leq^*\! \ul t} \Phi_\bb(\ul s) &\ \Rightarrow &\text{($\ul s$ closed monotone)} \\
      \EHAwst \vdash \Forall{\ul x} \Existsst{\ul y \leq^*\! \ul t} \Phi^\bb &\ \Rightarrow &\text{(characterisation)} \\
      \EHAwst + \PP \vdash \Forall{\ul x} \Existsst{\ul y \leq^*\! \ul t} \Phi &. 
    \end{align*}

    \item We have
    \begin{align*}
      \EHAwst + \PP \vdash \Phi_\mnexistsst \to \Existsst{\ul x} \Psi &\ \Rightarrow &\text{(soundness)} \\
      \EHAwst \vdash (\Phi_\mnexistsst \to \Existsst{\ul x} \Psi)_\bb(\ul t,\ul s) &\ \Rightarrow &\text{(lemma~\ref{lemma realisability of defined quantifications})} \\
      \EHAwst \vdash \Phi_\mnexistsst \to \Exists{\ul x \leq^*\! \ul t} \Psi_\bb(\ul s) &\ \Rightarrow &\text{($\ul t$ closed)} \\
      \EHAwst \vdash \Phi_\mnexistsst \to \Existsst{\ul x \leq^*\! \ul t} \Psi_\bb(\ul s) &\ \Rightarrow &\text{($\ul s$ closed monotone)} \\
      \EHAwst \vdash \Phi_\mnexistsst \to \Existsst{\ul x \leq^*\! \ul t} \Psi^\bb &\ \Rightarrow &\text{(characterisation)} \\
      \EHAwst + \PP \vdash \Phi_\mnexistsst \to \Existsst{\ul x \leq^*\! \ul t} \Psi &.
    \end{align*}

    \item Follows from point~\ref{realisability application mACw}.

    \item It is well known that there is an internal atomic formula~$\phi(x,y,z)$ of $\EHAw$ such that for all $x,y,z \in \mathbb N$ we have the following equivalence (where $\bar n$ denotes the numeral associated to $n \in \mathbb N$): $\EHAw \vdash \phi(\bar x,\bar y,\bar z)$ if and only if the Turing machine coded by $x$ when given input coded by $y$ halts with computation history coded by $z$. Let
    \begin{equation*}
      \Phi \defequiv (\Forallst z \neg\phi \to 0 =_0 1) \to \Existsst{\tilde z} (\Forall{z \leq_0 \tilde z} \neg\phi \to 0 =_0 1).
    \end{equation*}
    \begin{description}
      \item[$\EHAwst + \PP \nvdash \Phi$.] We have
      \begin{align*}
        \EHAwst + \PP \vdash \Phi &\ \Rightarrow &\text{(logic)} \\
        \EHAwst + \PP \vdash \Forall{x,y} \Phi &\ \Rightarrow &\text{(soundness)} \\
        \EHAwst \vdash (\Forallst{x,y} \Phi)_\bb(t) &\ \Rightarrow &\text{(lemma~\ref{lemma realisability of defined quantifications})} \\
        \EHAwst \vdash \mForall{\tilde x,\tilde y} \Forall{x,y \leq_0 \tilde x,\tilde y} (\neg\Forallst z \neg\phi \to \\
        \Exists{\tilde z \leq_0 t\tilde x\tilde y} \neg \Forall{z \leq_0 \tilde z} \neg\phi) &\ \Rightarrow &\text{(logic, arithmetic)} \\
        \EHAwst \vdash \Forallst{x,y} (\neg\Forallst z \neg\phi \to \neg\Forallst{z \leq_0 txy} \neg\phi) &.
      \end{align*}
      If $\EHAwst + \PP \vdash \Phi$, then there is a term~$t$ as above and $t$ induces a computable function that can be used to solve the halting problem, a contradiction.
      \item[$\EHAwst + \PP \nvdash \neg \Phi$.] If $\EHAwst + \PP \vdash \neg \Phi$, then $\EHAwst + \PP$ is inconsistent because $\EHAwst \vdash \neg\neg\Phi$.
    \end{description}

    \item\label{item LEM realisability applications} Let $\phi$ be as above and
    \begin{gather*}
      \Phi \defequiv \Existsst z \phi(x,y,z) \vee \neg\Existsst z \phi(x,y,z), \\
      \tilde \Phi \defequiv \Existsst z \phi(x,y,z) \vee \Forallst z \neg\phi(x,y,z).
    \end{gather*}
    \begin{description}
      \item[$\EHAwst + \PP \nvdash \Phi$.] We have
      \begin{align*}
        \EHAwst + \PP \vdash \Phi &\ \Rightarrow &\text{(logic)} \\
        \EHAwst + \PP \vdash \Forall{x,y} \tilde\Phi &\ \Rightarrow &\text{(soundness)} \\
        \EHAwst \vdash (\Forallst{x,y} \tilde\Phi)_\bb(t) &\ \Rightarrow &\text{(lemma~\ref{lemma realisability of defined quantifications})} \\
        \EHAwst \vdash \mForall{\tilde x,\tilde y} \forall{x,y \leq_0 \tilde x,\tilde y} \\
        (\Exists{z \leq_0 t\tilde x\tilde y} \phi(x,y,z) \vee \Forall z \neg\phi(x,y,z)) &\ \Rightarrow &\text{(logic, arithmetic)} \\
        \EHAwst \vdash \forallst{x,y} \\
        (\Existsst{z \leq_0 txy} \phi(x,y,z) \vee \Forallst z \neg\phi(x,y,z)) &.
      \end{align*}
      If $\EHAwst + \PP \vdash \Phi$, then there is a term~$t$ as above and $t$ induces a computable function that can be used to solve the halting problem, a contradiction.
      \item[$\EHAwst + \PP \nvdash \neg\Phi$.] If $\EHAwst + \PP \vdash \neg \Phi$, then $\EHAwst + \PP$ is inconsistent because $\EHAwst \vdash \neg\neg\Phi$.\qedhere
    \end{description}
  \end{enumerate}
\end{proof}

\section{Intuitionistic nonstandard bounded functional interpretation}

\begin{definition}
  The \emph{intuitionistic nonstandard bounded functional interpretation}~$\B$ assigns to each formula~$\Phi$ of $\EHAwst$ the formulas~$\Phi^\B$ and $\Phi_\B(\ul a;\ul b)$ of $\EHAwst$ such that $\Phi^\B \equiv \mExistsst{\ul a} \mForallst{\ul b} \Phi_\B(\ul a;\ul b)$ accordingly to the following clauses where $\Phi_\B(\ul a;\ul b)$ is the part inside square brackets. For atomic formulas, we define:
  \begin{enumerate}
    \item $\Phi^\B \defequiv [\Phi]$ for internal atomic formulas~$\Phi$ (so the tuples $\ul a$ and $\ul b$ are empty);
    \item $\st(t)^\B \defequiv \mExistsst a [t \leq^*\! a]$ (so the tuple $\ul b$ is empty).
  \end{enumerate}
  For the remaining formulas, if $\Phi^\B \equiv \mExistsst{\ul a} \mForallst{\ul b} \Phi_\B(\ul a;\ul b)$ and $\Psi^\B \equiv \mExistsst{\ul c} \mForallst{\ul d} \Psi_\B(\ul c;\ul d)$ then we define:
  \begin{enumerate}
    \setcounter{enumi}{2}
    \item $(\Phi \wedge \Psi)^\B \defequiv \mExistsst{\ul a,\ul c} \mForallst{\ul b,\ul d} [\Phi_\B(\ul a;\ul b) \wedge \Psi_\B(\ul c;\ul d)]$;
    \item $(\Phi \vee \Psi)^\B \defequiv \mExistsst{\ul a,\ul c} \mForallst{\ul e,\ul f} [\mForall{\ul b \leq^*\! \ul e} \Phi_\B(\ul a;\ul b) \vee \mForall{\ul d \leq^*\! \ul f} \Psi_\B(\ul c;\ul d)]$;
    \item $(\Phi \rightarrow \Psi)^\B \defequiv \mExistsst{\ul C,\ul B} \mForallst{\ul a,\ul d} [\mForall{\ul b \leq^*\! \ul B\ul a\ul d} \Phi_\B(\ul a;\ul b) \rightarrow \Psi_\B(\ul C\ul a;\ul d)]$;
    \item $(\Forall x \Phi)^\B \defequiv \mExistsst{\ul a} \mForallst{\ul b} [\Forall x \Phi_\B(\ul a;\ul b)]$;
    \item $(\Exists x \Phi)^\B \defequiv \mExistsst{\ul a} \mForallst{\ul c} [\Exists x \mForall{\ul b \leq^*\! \ul c} \Phi_\B(\ul a;\ul b)]$.
  \end{enumerate}
\end{definition}

We need the next lemma to prove the soundness and characterisation theorems of $\B$.

\begin{lemma}\label{lemma interpretation of defined quantifications}\mbox{}
  \begin{enumerate}
    \item For all internal formulas~$\phi$ of $\EHAwst$, we have $\phi^\B \equiv \phi_\B(\ul a) \equiv \phi$ (so the tuple $\ul a$ is empty).
    \item For all formulas~$\Phi$ of $\EHAwst$, we have (the equivalences being provable in $\EHAwst$):
    \begin{enumerate}
     \item $(\mForall{\ul x} \Phi)^\B \equiv \mExistsst{\ul A} \mForallst{\ul c,\ul b} [\mForall{\ul x \leq^*\! \ul c} \Phi_B(\ul A \ul c;\ul b)]$;
      \item $(\Forallst{\ul x} \Phi)^\B \equiv \mExistsst{\ul A} \mForallst{\ul c,\ul b} [\Forall{\ul x \leq^*\! \ul c} \Phi_B(\ul A \ul c;\ul b)]$;
      \item $(\mForallst{\ul x} \Phi)^\B \equiv \mExistsst{\ul A} \mForallst{\ul c,\ul b} [\mForall{\ul x \leq^*\! \ul c} \Phi_B(\ul A\ul c;\ul b)]$;
      \item $(\Forall{\ul x \leq^*\! \ul t} \Phi)^\B \equiv \mExistsst{\ul a} \mForallst{\ul b} [\Forall{\ul x \leq^*\! \ul t} \Phi_\B(\ul a;\ul b)]$;
      \item $(\mForall{\ul x \leq^*\! \ul t} \Phi)^\B \equiv \mExistsst{\ul a} \mForallst{\ul b} [\mForall{\ul x \leq^*\! \ul t} \Phi_\B(\ul a;\ul b)]$;
      \item $(\Existsst{\ul x} \Phi)^\B \leftrightarrow \mExistsst{\ul c,\ul a} \mForallst{\ul d} [\Exists{\ul x \leq^*\! \ul c} \mForall{\ul b \leq^*\! \ul d} \Phi_\B(\ul a;\ul b)]$;
      \item\label{item interpretation of defined quantifications} $(\mExistsst{\ul x} \Phi)^\B \leftrightarrow \mExistsst{\ul c,\ul a} \mForallst{\ul d} [\mExists{\ul x \leq^*\! \ul c} \mForall{\ul b \leq^*\! \ul d} \Phi_\B(\ul a;\ul b)]$;
      \item $(\Exists{\ul x \leq^*\! \ul t} \Phi)^\B \leftrightarrow \mExistsst{\ul a} \mForallst{\ul c} [\Exists{\ul x \leq^*\! \ul t} \mForall{\ul b \leq^*\! \ul c} \Phi_\B(\ul a;\ul b)]$;
      \item $(\mExists{\ul x \leq^*\! \ul t} \Phi)^\B \leftrightarrow \mExistsst{\ul a} \mForallst{\ul c} [\mExists{\ul x \leq^*\! \ul t} \mForall{\ul b \leq^*\! \ul c} \Phi_\B(\ul a;\ul b)]$.
    \end{enumerate}
  \end{enumerate}
\end{lemma}

\begin{proof}\mbox{}
  \begin{enumerate}
    \item The proof is by a simple induction on the length of $\phi$ (see the proof of lemma~\ref{lemma interpretation of internal formula} below).
    \item The proof is by simple calculations using that $\ul x \leq^*\! \ul t$ is an internal formula. For example, for point~\ref{item interpretation of defined quantifications}:
    \begin{align*}
      \Phi^\B &\equiv \mExistsst{\ul a} \mForallst{\ul b} \Phi_\bb(\ul a;\ul b), &\text{(assumption)} \\
      (\ul x \leq^*\! \ul x)^\B &\equiv \ul x \leq^*\! \ul x, &\text{($\ul x \leq^*\! \ul x$ internal)} \\
      \st(\ul x)^\B &\equiv \mExistsst{\ul c} [\ul x \leq^*\! \ul c], &\text{(calculation)} \\
      (\ul x \leq^*\! \ul x \wedge \st(\ul x))^\bb &\equiv \mExistsst{\ul c} [\ul x \leq^*\! \ul x \wedge \ul x \leq^*\! \ul c], &\text{(calculation)} \\
      (\ul x \leq^*\! \ul x \wedge \st(\ul x) \wedge \Phi)^\B &\equiv \mExistsst{\ul c,\ul a} \mForallst{\ul b} [\ul x \leq^*\! \ul x \wedge {} \\ &\phantom{{}\equiv{}} \ul x \leq^*\! \ul c \wedge \Phi_\B(\ul a;\ul b)], &\text{(calculation)} \\
     (\Exists{\ul x} (\ul x \leq^*\! \ul x \wedge \st(\ul x) \wedge \Phi))^\B &\equiv \mExistsst{\ul c,\ul a} \mForallst{\ul d} [\Exists{\ul x} \mforall{\ul b \leq^*\! \ul d} \\ &\phantom{{}\equiv{}} \ul x \leq^*\! \ul x \wedge \ul x \leq^*\! \ul c \wedge \Phi_\B(\ul a;\ul b)], &\text{(calculation)} \\
      (\Exists{\ul x} (\ul x \leq^*\! \ul x \wedge \st(\ul x) \wedge \Phi))^\B &\leftrightarrow \mExistsst{\ul c,\ul a} \mForallst{\ul d} [\mexists{\ul x \leq^*\! \ul c} \\ &\phantom{{}\leftrightarrow{}} \mForallst{\ul b \leq^*\! \ul d} \Phi_\B(\ul a;\ul b)], &\text{(logic)}
    \end{align*}
    where $\mExistsst{\ul x} \Phi \equiv \Exists{\ul x} (\ul x \leq^*\! \ul x \wedge \st(\ul x) \wedge \Phi)$.\qedhere
  \end{enumerate}
\end{proof}

The formulas $\Phi_\B(\ul a;\ul b)$ are monotone on $\ul a$ in the sense of the following lemma.

\begin{lemma}[monotonicity of $\B$]
  For all formulas $\Phi$ of $\EHAwst$, we have $\EHAwst \vdash \Phi_\B(\ul a;\ul b) \wedge \ul a \leq^*\! \ul{\tilde a} \rightarrow \Phi_\B(\ul{\tilde a};\ul b)$.
\end{lemma}

\begin{proof}
  The proof is by a simple induction on the length of $\Phi$. For example, in the case of $\rightarrow$, we have $(\Phi \rightarrow \Psi)_\B(\ul C,\ul B;\ul a,\ul d) \equiv \mForall{\ul b \leq^*\! \ul B\ul a\ul d} \Phi_\B(\ul a;\ul b) \rightarrow \Psi_\B(\ul C\ul a;\ul d)$, we have $\EHAwst \vdash \Psi_\B(\ul c;\ul d) \wedge \ul c \leq^*\! \ul{\tilde c} \rightarrow \Psi_\B(\ul{\tilde c};\ul d)$ by induction hypothesis, and we prove $(\mForall{\ul b \leq^*\! \ul B\ul a\ul d} \Phi_\B(\ul a;\ul b) \rightarrow \Psi_\B(\ul C\ul a;\ul d)) \wedge \ul C,\ul B \leq^*\! \ul{\tilde C},\ul{\tilde D} \rightarrow (\mForall{\ul b \leq^*\! \ul{\tilde B}\ul a\ul d} \Phi_\B(\ul a;\ul b) \rightarrow \Psi_\B(\ul{\tilde C}\ul a;\ul d))$ by noticing that $\ul C,\ul D \leq^*\! \ul{\tilde C},\ul{\tilde D}$ implies $\ul C \ul a,\ul B\ul a\ul d \leq^*\! \ul{\tilde C}\ul a,\ul{\tilde B}\ul a\ul d$ (for monotone $\ul a,\ul d$) and so the induction hypothesis implies $\Psi_\B(\ul C\ul a;\ul d) \rightarrow \Psi_\B(\ul{\tilde C}\ul a;\ul d)$ and trivially $\mForall{\ul b \leq^*\! \ul{\tilde B}\ul a\ul d} \Phi_\B(\ul a;\ul b) \to \mForall{\ul b \leq^*\! \ul B\ul a\ul d} \Phi_\B(\ul a;\ul b)$ (for monotone $\ul a,\ul d$).
\end{proof}

\begin{lemma}
  \label{lemma interpretation of internal formula}
  For all internal formulas~$\phi$ of $\EHAwst$, we have $\phi^\B \equiv \phi_\B(\ul a;\ul b)$ (so the tuples~$\ul a$ and $\ul b$ are empty) and $\phi_\B \equiv \phi$.
\end{lemma}

\begin{proof}
  The proof is by a simple induction on the length of $\phi$. For example, in the case of $\to$, we have $(\phi \rightarrow \psi)^\B \equiv \mExistsst{\ul C,\ul B} \mForallst{\ul a,\ul d} [\mForall{\ul b \leq^*\! \ul B\ul a\ul d} \phi_\B(\ul a;\ul b) \rightarrow \psi_\B(\ul C\ul a;\ul d)]$ (where $\ul C,\ul B$ and $\ul a,\ul d$ are empty by induction hypothesis), which is $\phi_\B \to \psi_\B$ (where $\phi_\B \equiv \phi$ and $\psi_\B \equiv \psi$ by induction hypothesis), which is $\phi \to \psi$.
\end{proof}

\begin{lemma}
  \label{lemma interpretation is internal formula}
  For all formulas~$\Phi$ of $\EHAwst$, the formula $\Phi_\B(\ul a;\ul b)$ is internal.
\end{lemma}

\begin{proof}
  The proof is by a simple induction on the length of $\Phi$. For example, in the case of $\st(t)$, we have $\st(t)_\B(a) \equiv t \leq^*\! a$, which is internal because it is built from formulas of the form $r \leq_0 s$ (which is an abbreviation for an appropriate internal atomic formula) by means of $\forall$, $\to$ and $\wedge$.
\end{proof}

The following seven principles, inspired by or even taken from Ferreira, Oliva and Gaspar's articles~\cite{FerreiraGaspar2015,FerreiraOliva2005}, play an important role in the sequel as they are the characteristic principles of $\B$.

\begin{definition}
  We define the following principles:
  \begin{enumerate}
    \item the \emph{monotone choice}~$\mACw$ is all instances of $\mForallst x \mExistsst y \Phi \rightarrow \mExistsst Y \mForallst x \mExists{y \leq^*\! Yx} \Phi$;
    \item the \emph{realisation}~$\Rw$ is all instances of $\Forall x \Existsst y \Phi \rightarrow \mExistsst z \Forall x \Exists{y \leq^*\! z} \Phi$;
    \item the \emph{idealisation}~$\Iw$ is all instances of $\mForallst{\ul z} \Exists{\ul x} \Forall{\ul y \leq^*\! \ul z} \phi \rightarrow \Exists{\ul x} \Forallst{\ul y} \phi$;
    \item the \emph{independence of premises}~$\IPmforallwst$ is all instances of $(\mForallst x \phi \rightarrow \mExistsst y \Psi) \rightarrow \mExistsst z (\mForallst x \phi \rightarrow \mExists{y \leq^*\! z} \Psi)$;
    \item \emph{Markov's principle}~$\Mw$ is all instances of $(\mForallst x \phi \rightarrow \psi) \rightarrow \mExistsst y (\mForall{x \leq^*\! y} \phi \rightarrow \psi)$;
    \item the \emph{bounded universal disjunction principle}~$\BUDw$ are all instances of $\mForallst{\ul u,\ul v} (\Forall{\ul x \leq^*\! \ul u} \phi \vee \Forall{\ul y \leq^*\! \ul v} \psi) \rightarrow \Forallst{\ul x} \phi \vee \Forallst{\ul y} \psi$;
    \item the \emph{majorisability axioms}~$\MAJw$ are all instances of $\Forallst x \Existsst y (x \leq^*\! y)$.
  \end{enumerate}
\end{definition}

The principles above, where there are single variables $x$, $y$, $Y$ and $z$ instead of tuples of variables $\ul x$, $\ul y$, $\ul Y$ and $\ul z$, generalise to tuples of variables, that is each principle provably implies in $\EHAwst$ its variant where the single variables~$x$, $y$, $Y$ and $z$ are replaced by (respectively) the tuples of variables~$\ul x$, $\ul y$, $\ul Y$ and $\ul z$ (proof sketch: first we generalise from $x$ to $\ul x$ by induction on the length of $\ul x$ and then we generalise from $y$, $Y$ and $z$ to $\ul y$, $\ul Y$ and $\ul z$ by induction on the length of $\ul y$, $\ul Y$ and $\ul z$). The principles above, where there are tuples of variables $\ul u$, $\ul v$, $\ul x$, $\ul y$ and $\ul z$ instead of single variables $u$, $v$, $x$, $y$ and $z$, seemingly do \emph{not} generalise to tuples of variables (by induction on the length of the tuples).

The principle~$\Iw$ implies the principle $\BUDw$, that is $\EHAwst + \Iw$ proves all instances of $\BUDw$ (proof sketch: rewrite the disjunction in the premises of $\BUDw$ as an existential conjunction of implications getting $\mForallst{\ul u,\ul v} \Exists x \big((x =_0 0 \to \Forall{\ul x \leq^*\! \ul u} \phi) \wedge (x \neq_0 0 \to \Forall{\ul y \leq^*\! \ul v} \psi))\big)$, that is $\mForallst{\ul u,\ul v} \Exists x \Forall{\ul x,\ul y \leq^*\! \ul u,\ul v} ((x =_0 0 \to \phi) \wedge (x \neq_0 0 \to \psi))$, apply $\Iw$ getting $\Exists x \Forallst{\ul x,\ul y} ((x =_0 0 \to \phi) \wedge (x \neq_0 0 \to \psi))$, that is $\Exists x ((x =_0 0 \to \Forallst{\ul x} \phi) \wedge (x \neq_0 0 \to \Forallst{\ul y} \psi))$, and rewrite the existential conjunction of implications back as a disjunction getting the conclusion of $\BUDw$). The principle~$\Rw$ implies the principle $\MAJw$ (see page~\ref{page Rw implies MAJw}). However, we keep explicitly mentioning $\BUDw$ and $\MAJw$ (even in the presence of $\Iw$ and $\Rw$) because of their importance.

Now we present our main theorem about~$\B$. This theorem shows that $\B$:
\begin{enumerate}
  \item interprets $\EHAwst + \mACw + \Rw + \Iw + \IPmforallwst + \Mw + \BUDw + \MAJw$ into $\EHAwst$;
  \item extracts computational information~$\ul t$ from proofs in $\EHAwst + \mACw + \Rw + \Iw + \IPmforallwst + \Mw + \BUDw + \MAJw$.
\end{enumerate}

\begin{theorem}[soundness theorem of $\B$]
  For all formulas $\Phi$ of $\EHAwst$, if
  \begin{equation*}
     \EHAwst + \mACw + \Rw + \Iw + \IPmforallwst + \Mw + \BUDw + \MAJw \vdash \Phi,
  \end{equation*}
  then there are \emph{closed} monotone terms $\ul t$ of appropriate types such that
  \begin{equation*}
    \EHAwst \vdash \mForallst{\ul b} \Phi_\B(\ul t;\ul b).
  \end{equation*}
\end{theorem}

\begin{proof}
  The proof is by induction on length of the derivation of $\Phi$. The logical and arithmetical axioms and rules are dealt with in a way similar to the bounded functional interpretation~\cite{FerreiraOliva2005}. So let us focus on the standardness axioms and the characteristic principles (except $\BUDw$ and $\MAJw$, which follow respectively from $\Iw$ and $\Rw$). We will sometimes use implicitly lemmas~\ref{lemma interpretation of defined quantifications} and \ref{lemma interpretation of internal formula}.
  \begin{description}
    \item[$x = y \wedge \st(x) \rightarrow \st(y)$.] Its interpretation is $\mExistsst B \mForallst a (x = y \wedge x \leq^*\! a \to y \leq^*\! Ba)$. It is interpreted by $B \defeq \la{a}{a}$.

    \item[$\st(y) \wedge x \leq^*\! y \rightarrow \st(x)$.] Its interpretation is $\mExistsst B \mForallst a (y \leq^*\! a \wedge x \leq^*\! y \to x \leq^*\! Ba)$. It is interpreted by $B \defeq \la{a}{a}$.

    \item[$\st(t)$.] Its interpretation is $\mExistsst a (t \leq^*\! a)$. It is interpreted by a closed term $\tilde t$ such that $t \leq^*\! \tilde t$, which exists by Howard's majorisability theorem.

    \item[$\st(x) \wedge \st(y) \rightarrow \st(xy)$.] Its interpretation is $\mExistsst C \mForallst {a,b} (x \leq^*\! a \wedge y \leq^*\! b \to xy \leq^*\! Cab)$. It is interpreted by $C \defeq \la{a,b}{ab}$.

    \item[$\mForallst x \mExistsst y \Phi \rightarrow \mExistsst Y \mForallst x \mExists{y \leq^*\! Yx} \Phi$.] Its interpretation is equivalent to
    \begin{equation*}
      \begin{gathered}
        \mExistsst{J,\ul{\mathrm C},G,\ul F} \mforallst{E,\ul A,k,\ul l} \\
        (\mForall{g,\ul f \leq^*\! GE\ul A k\ul l,\ul F E\ul A k\ul l} \mForall{x \leq^*\! g} \mExists{y \leq^*\! Eg} \mForall{\ul b \leq^*\! \ul f} \Phi_B(\ul A g;\ul b) \to \\
        \Exists{Y \leq^*\! JE\ul A} \mForall{i,\ul h \leq^*\! k,\ul l} \mForall{x \leq^*\! i} \mExists{y \leq^*\! Yx} \mForall{\ul d \leq^*\! \ul h} \Phi_\B(\ul C E\ul A i; \ul d)).
      \end{gathered}
    \end{equation*}
    It is interpreted by $J \defeq \la{E,\ul A}{E}$, $\ul{\mathrm C} \defeq \la{E,\ul A,i}{\ul A i}$, $G \defeq \la{E,\ul A,k,\ul l}{k}$ and $\ul F \defeq \la{E,\ul A,k,\ul l}{\ul l}$. Indeed, from the premise~$\mForall{g,\ul f \leq^*\! k,\ul l} \mForall{x \leq^*\! g} \mExists{y \leq^*\! Eg} \mForall{\ul b \leq^*\! \ul f} \Phi_B(\ul A g;\ul b)$ with $x \defeq g$ we get $\mForall{x,\ul f \leq^*\! k,\ul l} \mExists{y \leq^*\! Ex} \mForall{\ul b \leq^*\! \ul f} \Phi_B(\ul A g;\ul b)$, so we get the conclusion~$\Exists{Y \leq^*\! E} \mForall{i,\ul h \leq^*\! k,\ul l} \mForall{x \leq^*\! i} \mExists{y \leq^*\! Yx} \mForall{\ul d \leq^*\! \ul h} \Phi_\B(\ul A i; \ul d))$ with $Y \defeq E$.

    \item[$\Forall x \Existsst y \Phi \rightarrow \mExistsst z\Forall x \Exists{y \leq^*\! z} \Phi$.] Its interpretation is equivalent to
    \begin{equation*}
      \begin{gathered}
        \mExistsst{H,\ul C,\ul F} \mforallst{e,\ul a,\ul i} \\
        (\mForall{\ul f \leq^*\! \ul F e\ul a i} \Forall x \Exists{y \leq^*\! e} \mForall{\ul b \leq^*\! \ul f} \Phi_\B(\ul a;\ul b) \to \\
        \mExists{z \leq^*\! He\ul a} \mForall{\ul g \leq^*\! \ul i} \Forall x \Exists{y \leq^*\! z} \mForall{\ul d \leq^*\! \ul g} \Phi_\B(\ul C e\ul a;\ul d)).
      \end{gathered}
    \end{equation*}
    It is interpreted by $H \defeq \la{e,\ul a}{e}$, $\ul C \defeq \la{e,\ul a}{\ul a}$ and $\ul F \defeq \la{e,\ul a,\ul i}{\ul i}$.

    \item[$\mForallst{\ul z} \Exists{\ul x} \Forall{\ul y \leq^*\! \ul z} \phi \rightarrow \Exists{\ul x} \Forallst{\ul y} \phi$.] Its interpretation is
    \begin{equation*}
      \mExistsst{\ul A} \mForallst{\ul c} (\mForall{\ul a \leq^*\! \ul A\ul c} \mForall{\ul z \leq^*\! \ul a} \Exists{\ul x} \Forall{\ul y \leq^*\! \ul c} \phi \to \Exists{\ul x} \mForall{\ul b \leq^*\! \ul c} \mForall{\ul y \leq^*\! \ul b} \phi).
    \end{equation*}
    It is interpreted by $\ul A \defeq \la{\ul c}{\ul c}$.

  \item[$(\mForallst x \phi \rightarrow \mExistsst y \Psi) \rightarrow \mExistsst z (\mForallst x \phi \rightarrow \mExists{y \leq^*\! z} \Psi)$.] Its interpretation is equivalent to
    \begin{equation*}
      \begin{gathered}
        \mExistsst{J,\ul C,\ul{\mathrm H},\ul G} \mForallst{f,\ul a,\ul E,\ul k} \\
        \big(\mForall{\ul g \leq^*\! \ul G f\ul a\ul E\ul k} (\mForall{\ul e \leq^*\! \ul E\ul g} \mForall{\ul x \leq^*\! \ul e} \phi \to \mExists{y \leq^*\! f} \mForall{\ul b \leq^*\! \ul g} \Psi_\B(\ul a;\ul b)) \\
        \mExists{z \leq^*\! Jf\ul a\ul E} \mForall{\ul i \leq^*\! \ul k} (\mForall{\ul h \leq^*\! \ul{\mathrm H}f\ul a\ul i} \mForall{\ul x \leq^*\! \ul h} \phi \to \mExists{y \leq^*\! z} \mForall{\ul d \leq^*\! \ul i} \Psi_\B(\ul Cf\ul a\ul E;\ul d))\big).
      \end{gathered}
    \end{equation*}
    It is interpreted by $J \defeq \la{f,\ul a,\ul E}{f}$, $\ul C \defeq \la{f,\ul a,\ul E}{\ul a}$, $\ul{\mathrm H} \defeq \la{f,\ul a,\ul E,\ul i}{\ul E\ul i}$ and $\ul G \defeq \la{f,\ul a,\ul E,\ul k}{\ul k}$.

   \item[$(\mForallst x \phi \rightarrow \psi) \rightarrow \mExistsst y (\mForall{x \leq^*\! y} \phi \rightarrow \psi)$.] Its interpretation is equivalent to
    \begin{equation*}
      \mExistsst{\ul C} \mForallst{\ul b} ((\mForall{\ul a \leq^*\! \ul b} \mForall{\ul x \leq^*\! \ul a} \phi \to \psi) \to (\mExists{\ul y \leq^*\! \ul C\ul b} (\mForall{\ul x \leq^*\! \ul y} \phi \to \psi)).
    \end{equation*}
    It is interpreted by $\ul C \defeq \la{\ul b}{\ul b}$.\qedhere
  \end{description}
\end{proof}

\begin{theorem}[characterisation theorem of $\B$] For all formulas $\Phi$ of $\EHAwst$, we have
  \begin{equation*}
    \EHAwst + \mACw + \Rw + \Iw + \IPmforallwst + \Mw + \BUDw + \MAJw \vdash \Phi \leftrightarrow \Phi^\B.
  \end{equation*}
\end{theorem}

\begin{proof}
  The proof is by induction on the length of $\Phi$. The cases of internal atomic formulas and $\wedge$ are easy, so let us focus on the remaining cases. We will sometimes use implicitly lemma~\ref{lemma interpretation is internal formula}.
  \begin{description}
    \item[$\st(t)$.] We have
    \begin{align*}
      \st(t) &\leftrightarrow &\text{($\MAJw$, $\st(y) \wedge x \leq^*\! y \rightarrow \st(x)$)} \\
      \mExistsst a (t \leq^*\! a) &\equiv &\text{(definition)} \\
      \st(t)^\B &.
    \end{align*}

    \item[$\vee$.] We have
    \begin{align*}
      \Phi \vee \Psi &\leftrightarrow &\text{(induction hypothesis)} \\
      \mExistsst{\ul a} \mForallst{\ul b} \Phi_\B(\ul a;\ul b) \vee \mExistsst{\ul c} \mForallst{\ul d} \Psi_\B(\ul c;\ul d) &\leftrightarrow &\text{(logic)} \\
      \mExistsst{\ul a,\ul c} (\mForallst{\ul b} \Phi_\B(\ul a;\ul b) \vee \mForallst{\ul d} \Psi_\B(\ul c;\ul d)) &\leftrightarrow &\text{($\BUDw$)} \\
      \mExistsst{\ul a,\ul c} \mForallst{\ul e,\ul f} (\mForall{\ul b \leq^*\! \ul e} \Phi_\B(\ul a;\ul b) \vee \mForall{\ul d \leq^*\! \ul f} \Psi_\B(\ul c;\ul d)) &\equiv &\text{(definition)} \\
      (\Phi \vee \Psi)^\B &.
    \end{align*}

    \item[$\rightarrow$.] We have
    \begin{align*}
      (\Phi \rightarrow \Psi) &\leftrightarrow &\text{(induction hypothesis)} \\
      (\mExistsst{\ul a} \mForallst{\ul b} \Phi_\B(\ul a;\ul b) \to \mExistsst{\ul c} \mForallst{\ul d} \Psi_\B(\ul c;\ul d)) &\leftrightarrow &\text{(logic)} \\
      \mForallst{\ul a} (\mForallst{\ul b} \Phi_\B(\ul a;\ul b) \to \mExistsst{\ul c} \mForallst{\ul d} \Psi_\B(\ul c;\ul d)) &\leftrightarrow &\text{($\IPmforallwst$)} \\
\mForallst{\ul a} \mExistsst{\ul c} (\mForallst{\ul b} \Phi_\B(\ul a;\ul b) \to \mExistsst{\ul{\tilde c} \leq^*\! \ul c} \mForallst{\ul d} \Psi_\B(\ul{\tilde c};\ul d)) &\leftrightarrow &\text{(monotonicity)} \\
\mForallst{\ul a} \mExistsst{\ul c} (\mForallst{\ul b} \Phi_\B(\ul a;\ul b) \to \mForallst{\ul d} \Psi_\B(\ul c;\ul d)) &\leftrightarrow &\text{(logic)} \\
    \mForallst{\ul a} \mExistsst{\ul c} \mForallst{\ul d} (\mForallst{\ul b} \Phi_\B(\ul a;\ul b) \to \Psi_\B(\ul c;\ul d)) &\leftrightarrow &\text{($\Mw$)} \\
    \mForallst{\ul a} \mExistsst{\ul c} \mForallst{\ul d} \mExistsst{\ul{\tilde b}} (\mForallst{\ul{\tilde b} \leq^*\! \ul b} \Phi_\B(\ul a;\ul b) \to \Psi_\B(\ul c;\ul d)) &\leftrightarrow &\text{($\mACw$)} \\
    \mExistsst{\ul C,\ul B} \mForallst{\ul a} \mExistsst{\ul c \leq^*\! \ul C\ul a} \mForallst{\ul d} \mExistsst{\ul{\tilde b} \leq^*\! \ul B\ul a\ul d} \\ (\mForallst{\ul b \leq^*\! \ul{\tilde b}} \Phi_\B(\ul a;\ul b) \to \Psi_\B(\ul c;\ul d)) &\leftrightarrow &\text{(logic, monotonicity)} \\
    \mExistsst{\ul C,\ul B} \mForallst{\ul a,\ul d} (\mForallst{\ul b \leq^*\! B\ul a\ul d} \Phi_\B(\ul a;\ul b) \to \Psi_\B(\ul C\ul a;\ul d)) &\equiv &\text{(definition)} \\
      (\Phi \to \Psi)^\B &.
    \end{align*}

    \item[$\forall$.] We have
    \begin{align*}
      \Forall x \Phi &\leftrightarrow &\text{(induction hypothesis)} \\
      \Forall x \mExistsst{\ul a} \mForallst{\ul b} \Phi_\B(\ul a;\ul b) &\leftrightarrow &\text{($\Rw$)} \\
      \mExistsst{\ul a} \Forall x \mExists{\ul{\tilde a} \leq^*\! \ul a} \mForallst{\ul b} \Phi_\B(\ul{\tilde a};\ul b) &\leftrightarrow &\text{(monotonicity)} \\
      \mExistsst{\ul a} \Forall x \mForallst{\ul b} \Phi_\B(\ul a;\ul b) &\leftrightarrow &\text{(logic)} \\
      \mExistsst{\ul a} \mForallst{\ul b} \Forall x \Phi_\B(\ul a;\ul b) &\equiv &\text{(definition)} \\
      (\Forall x \Phi)^\B &.
    \end{align*}

    \item[$\exists$.] We have
    \begin{align*}
      \Exists x \Phi &\leftrightarrow &\text{(induction hypothesis)} \\
      \Exists x \mExistsst{\ul a} \mForallst{\ul b} \Phi_\B(\ul a;\ul b) &\leftrightarrow &\text{(logic)} \\
      \mExistsst{\ul a} \Exists x \mForallst{\ul b} \Phi_\B(\ul a;\ul b) &\leftrightarrow &\text{($\Iw$)} \\
\mExistsst{\ul a} \mForallst{\ul c} \Exists x \mForall{\ul b \leq^*\! \ul c} \Phi_\B(\ul a;\ul b) &\equiv &\text{(definition)} \\
      (\Exists x \Phi)^\B &. &&\qedhere
    \end{align*}
  \end{description}
\end{proof}

In the next result we give some applications of $\B$:
\begin{enumerate}
  \item the first, second, third, fourth, fifth, seventh and eight points are similar to the applications of $\bb$;
  \item the sixth point is a term-and-rule variant of a generalisation $(\mForallst x \phi \rightarrow \Existsst{\ul z} \psi) \rightarrow \mExistsst y (\mForall{x \leq^*\! y} \phi \rightarrow \Existsst{\ul z} \psi)$ of $\Mw$ (since $\Mw$ is the particular case in which the tuple~$\ul z$ is empty).
\end{enumerate}

\begin{application}
  Let $\EHAwst + \PP \defeq \EHAwst + \mACw + \Rw + \Iw + \IPmforallwst + \Mw + \BUDw + \MAJw$.
  \begin{enumerate}
    \item The theories $\EHAwst + \PP$ and $\EHAwst$ are equiconsistent.
    \item If $\EHAwst + \PP \vdash \Existsst{\ul x} \Phi$, then there are closed monotone terms~$\ul t$ of appropriate types such that $\EHAwst + \PP \vdash \Existsst{\ul x \leq^*\! \ul t} \Phi$. Moreover, if $\Phi$ is internal, then in the conclusion we can replace $\EHAwst + \PP$ by $\EHAwst$. Analogously, if all quantifications of $\ul x$ are monotone.
    \item\label{interpretation application mACw} If $\EHAwst + \PP \vdash \Forallst{\ul x} \Existsst{\ul y} \Phi$, then there are closed monotone terms~$\ul t$ of appropriate types such that $\EHAwst + \PP \vdash \mForallst{\ul z} \Forall{\ul x \leq^*\! \ul z} \Existsst{\ul y \leq^*\! \ul t\ul z} \Phi$. Also, in the conclusion we can replace $\mForallst{\ul z} \Forall{\ul x \leq^*\! \ul z} \Existsst{\ul y \leq^*\! \ul t\ul z} \Phi$ by $\mForallst{\ul x} \Existsst{\ul y \leq^*\! \ul t\ul x} \Phi$. Moreover, if $\Phi$ is internal, then in the conclusion we can replace $\EHAwst + \PP$ by $\EHAwst$. Analogously, if all quantifications of $\ul x$ are monotone, or if all quantifications of $\ul y$ are monotone, or if both.
    \item If $\EHAwst + \PP \vdash \Forall{\ul x} \Existsst{\ul y} \Phi$, then there are closed monotone terms~$\ul t$ of appropriate types such that $\EHAwst + \PP \vdash \Forall{\ul x} \Existsst{\ul y \leq^*\! \ul t} \Phi$. Moreover, if $\Phi$ is internal, then in the conclusion we can replace $\EHAwst + \PP$ by $\EHAwst$. Analogously, if all quantifications of $\ul x$ are monotone, or if all quantifications of $\ul y$ are monotone, or if both.
    \item If $\EHAwst + \PP \vdash \Forallst{\ul x} \phi \to \Existsst{\ul y} \Psi$, then there are closed monotone terms~$\ul t$ of appropriate types such that $\EHAwst + \PP \vdash \Forallst{\ul x} \phi \to \Existsst{\ul y \leq^*\! \ul t} \Psi$. Moreover, if $\Psi$ is internal, then in the conclusion we can replace $\EHAwst + \PP$ by $\EHAwst$. Analogously, if all quantifications of $\ul x$ are monotone, or if all quantifications of $\ul y$ are monotone, or if both.
    \item If $\EHAwst + \PP \vdash \Forallst{\ul x} \phi \to \Existsst{\ul y} \psi$, then there are closed monotone terms~$\ul s,\ul t$ of appropriate types such that $\EHAwst \vdash \Forallst{\ul x \leq^*\! \ul s} \phi \to \Existsst{\ul y \leq^*\! \ul t} \psi$. Analogously, if all quantifications of $\ul x$ are monotone, or if all quantifications of $\ul y$ are monotone, or if both.
    \item If $\EHAwst + \PP \vdash \mForallst{\ul x} \Existsst{\ul y} \phi$, then $\EHAwst \vdash  \mForallst{\ul x} \Existsst{\ul y} \phi$. Analogously, if all quantifications of $\ul y$ are monotone.
    \item If $\EHAwst + \PP$ is consistent, then there is an instance~$\Phi$ of $\LEM$ such that $\EHAwst + \PP \nvdash \Phi$ and $\EHAwst + \PP \nvdash \neg\Phi$.
  \end{enumerate}
\end{application}

\begin{proof}
  The ``Moreover, \ldots'' parts are proved analogously to their preceding parts but using lemma~\ref{lemma interpretation of internal formula} instead of the characterisation theorem to remain in $\EHAwst$ instead of going to $\EHAwst + \PP$. The ``Analogously, \ldots'' parts are proved analogously to their proceedings parts and sometimes are also a corollary to their preceding parts. We will sometimes use implicitly lemma~\ref{lemma interpretation is internal formula}. We will sometimes use the fact that $\Iw$ implies $\mForallst{\ul z} \Exists{\ul x \leq^*\! \ul t} \mForall{\ul y \leq^*\! \ul z} \tilde\phi \rightarrow \Exists{\ul x \leq^*\! \ul t} \mForallst{\ul y} \tilde\phi$ (proof sketch: rewrite the premise as $\mForallst{\ul z} \Exists{\ul x} \Forall{\ul y \leq^*\! \ul z} (\ul x \leq^*\! \ul t \wedge (\ul y \leq^*\! \ul y \to \tilde\phi))$, apply $\Iw$ with $\phi \defequiv \ul x \leq^*\! \ul t \wedge (\ul y \leq^*\! \ul y \to \tilde\phi)$ getting $\Exists{\ul x} \Forallst{\ul y} (\ul x \leq^*\! \ul t \wedge (\ul y \leq^*\! \ul y \to \tilde\phi))$, and then rewrite this as the conclusion). We will sometimes use implicitly that closed terms are standard.
  \begin{enumerate}
    \item We have
    \begin{align*}
      \EHAwst + \PP \vdash 0 =_0 1 &\ \Rightarrow & \text{(soundness)} \\
      \EHAwst \vdash (0 =_0 1)_\B &\ \Rightarrow & \text{(definition)} \\
      \EHAwst \vdash 0 =_0 1 &.
    \end{align*}
    Trivially, $\EHAwst \vdash 0 =_0 1 \ \Rightarrow \ \EHAwst + \PP \vdash 0 =_0 1$.

    \item We have
    \begin{align*}
      \EHAwst + \PP \vdash \Existsst{\ul x} \Phi &\ \Rightarrow &\text{(soundness)} \\
      \EHAwst \vdash (\Existsst{\ul x} \Phi)_\B(\ul t,\ul s) &\ \Rightarrow &\text{(lemma~\ref{lemma interpretation of defined quantifications})} \\
      \EHAwst \vdash \mForallst{\ul d} \Exists{\ul x \leq^*\! \ul t} \mForall{\ul b \leq^*\! \ul d} \Phi_\B(\ul s;\ul b) &\ \Rightarrow &\text{($\Iw$)} \\
      \EHAwst + \Iw \vdash \Exists{\ul x \leq^*\! \ul t} \mForall{\ul b} \Phi_\B(\ul s;\ul b) &\ \Rightarrow &\text{($\ul t$ closed)} \\
      \EHAwst + \Iw \vdash \Existsst{\ul x \leq^*\! \ul t} \mForall{\ul b} \Phi_\B(\ul s;\ul b) &\ \Rightarrow &\text{($\ul s$ closed monotone)} \\
      \EHAwst + \Iw \vdash \Existsst{\ul x \leq^*\! \ul t} \Phi^\B &\ \Rightarrow &\text{(characterisation)} \\
      \EHAwst + \PP \vdash \Existsst{\ul x \leq^*\! \ul t} \Phi &. \\
    \end{align*}

    \item We have
    \begin{align*}
      \EHAwst + \PP \vdash \Forallst{\ul x} \Existsst{\ul y} \Phi &\ \Rightarrow &\text{(soundness)} \\
      \EHAwst \vdash (\Forallst{\ul x} \Existsst{\ul y} \Phi)_\B(\ul t,\ul s) &\ \Rightarrow &\text{(lemma~\ref{lemma interpretation of defined quantifications})} \\
      \EHAwst \vdash \mForallst{\ul z,\ul d} \Forall{\ul x \leq^*\! \ul z} \Exists{\ul y \leq^*\! \ul t\ul z} \mForall{\ul b \leq^*\! \ul d} \Phi_\B(\ul s\ul z;\ul d) &\ \Rightarrow &\text{($\Iw$)} \\
      \EHAwst + \Iw \vdash \mForallst{\ul z} \Forall{\ul x \leq^*\! \ul z} \Exists{\ul y \leq^*\! \ul t\ul z} \mForallst{\ul b} \Phi_\B(\ul s\ul z;\ul d) &\ \Rightarrow &\text{($\ul t$ closed)} \\
      \EHAwst + \Iw \vdash \mForallst{\ul z} \Forall{\ul x \leq^*\! \ul z} \Existsst{\ul y \leq^*\! \ul t\ul z} \mForallst{\ul b} \Phi_\B(\ul s\ul z;\ul d) &\ \Rightarrow &\text{($\ul s\ul z$ closed monotone)} \\
      \EHAwst + \Iw \vdash \mForallst{\ul z} \Forall{\ul x \leq^*\! \ul z} \Existsst{\ul y \leq^*\! \ul t\ul z} \Phi^\B &\ \Rightarrow &\text{(characterisation)} \\
      \EHAwst + \PP \vdash \mForallst{\ul z} \Forall{\ul x \leq^*\! \ul z} \Existsst{\ul y \leq^*\! \ul t\ul z} \Phi &.
    \end{align*}
    Also, $\EHAwst + \PP \vdash \mForallst{\ul z} \Forall{\ul x \leq^*\! \ul z} \Existsst{\ul y \leq^*\! \ul t\ul z} \Phi \ \Rightarrow \ \EHAwst + \PP \vdash \mForallst{\ul x} \Existsst{\ul y \leq^*\! \ul t\ul x} \Phi$ (by taking $\ul x \defeq \ul z$).

    \item We have
    \begin{align*}
      \EHAwst + \PP \vdash \Forall{\ul x} \Existsst{\ul y} \Phi &\ \Rightarrow &\text{(soundness)} \\
      \EHAwst \vdash (\Forall{\ul x} \Existsst{\ul y} \Phi)_\B(\ul t,\ul s) &\ \Rightarrow &\text{(lemma~\ref{lemma interpretation of defined quantifications})} \\
      \EHAwst \vdash \mForallst{\ul d} \Forall{\ul x} \Exists{\ul y \leq^*\! \ul t} \mForall{\ul b \leq^*\! \ul d} \Phi_\B(\ul s;\ul b) &\ \Rightarrow &\text{($\Iw$)} \\
      \EHAwst + \Iw \vdash \Forall{\ul x} \Exists{\ul y \leq^*\! \ul t} \mForallst{\ul b} \Phi_\B(\ul s;\ul b) &\ \Rightarrow &\text{($\ul t$ closed)} \\
      \EHAwst + \Iw \vdash \Forall{\ul x} \Existsst{\ul y \leq^*\! \ul t} \mForallst{\ul b} \Phi_\B(\ul s;\ul b) &\ \Rightarrow &\text{($\ul s$ closed monotone)} \\
      \EHAwst + \Iw \vdash \Forall{\ul x} \Existsst{\ul y \leq^*\! \ul t} \Phi^\B &\ \Rightarrow &\text{(characterisation)} \\
      \EHAwst + \PP \vdash \Forall{\ul x} \Existsst{\ul y \leq^*\! \ul t} \Phi &.
    \end{align*}

    \item We have
    \begin{align*}
      \EHAwst + \PP \vdash \Forallst{\ul x} \phi \to \Existsst{\ul y} \Psi &\ \Rightarrow &\text{(soundness)} \\
      \EHAwst \vdash (\Forallst{\ul x} \phi \to \Existsst{\ul y} \Psi)_\B(\ul t,\ul r,\ul s) &\ \Rightarrow &\text{(lemma~\ref{lemma interpretation of defined quantifications})} \\
      \EHAwst \vdash \mForallst{\ul e} \\
      (\mForall{\ul c \leq^*\! \ul s\ul e} \Forall{\ul x \leq^*\! \ul c} \phi \to \Exists{\ul y \leq^*\! \ul t} \mForall{\ul b \leq^*\! \ul e} \Psi_\B(\ul r;\ul b)) &\ \Rightarrow &\text{($\ul s\ul e$ standard)} \\
      \EHAwst \vdash \Forallst{\ul x} \phi \to \mForallst{\ul e} \Exists{\ul y \leq^*\! \ul t} \mForall{\ul b \leq^*\! \ul e} \Psi_\B(\ul r;\ul b)) &\ \Rightarrow &\text{($\Iw$)} \\
      \EHAwst + \Iw \vdash \Forallst{\ul x} \phi \to \Exists{\ul y \leq^*\! \ul t} \mForallst{\ul b} \Psi_\B(\ul r;\ul b)) &\ \Rightarrow &\text{($\ul t$ closed)} \\
      \EHAwst + \Iw \vdash \Forallst{\ul x} \phi \to \Existsst{\ul y \leq^*\! \ul t} \mForallst{\ul b} \Psi_\B(\ul r;\ul b)) &\ \Rightarrow &\text{($\ul r$ closed monotone)} \\
      \EHAwst + \Iw \vdash \Forallst{\ul x} \phi \to \Existsst{\ul y \leq^*\! \ul t} \Psi^\B &\ \Rightarrow &\text{(characterisation)} \\
      \EHAwst + \PP \vdash \Forallst{\ul x} \phi \to \Existsst{\ul y \leq^*\! \ul t} \Psi &.
    \end{align*}

    \item We have
    \begin{align*}
      \EHAwst + \PP \vdash \Forallst{\ul x} \phi \to \Existsst{\ul y} \psi &\ \Rightarrow &\text{(soundness)} \\
      \EHAwst \vdash (\Forallst{\ul x} \phi \to \psi)_\B(\ul t,\ul s;\langle\rangle) &\ \Rightarrow &\text{(lemma~\ref{lemma interpretation of defined quantifications})} \\
      \EHAwst \vdash \mForall{\ul a \leq^*\! \ul s} \Forall{\ul x \leq^*\! \ul a} \phi \to \Exists{\ul y \leq^*\! \ul t} \psi &\ \Rightarrow &\text{($\ul s,\ul t$ closed monotone)} \\
      \EHAwst \vdash \Forallst{\ul x \leq^*\! \ul s} \phi \to \Existsst{\ul y \leq^*\! \ul t} \psi &.
    \end{align*}

    \item Follows from point~\ref{interpretation application mACw}.

    \item Let $\phi$, $\Phi$ and $\tilde\Phi$ be as in point~\ref{item LEM realisability applications} in the proof of application~\ref{realisability applications}.
    \begin{description}
      \item[$\EHAwst + \PP \nvdash \Phi$.] We have
      \begin{align*}
        \EHAwst + \PP \vdash \Phi &\ \Rightarrow &\text{(logic)} \\
        \EHAwst + \PP \vdash \Forall{x,y} \tilde\Phi &\ \Rightarrow &\text{(soundness)} \\
        \EHAwst \vdash \mForallst{\tilde x,\tilde y,\tilde z} (\Forallst{x,y} \tilde\Phi)_\B(t;\tilde x,\tilde y,\tilde z) &\ \Rightarrow &\text{(lemma~\ref{lemma realisability of defined quantifications})} \\
        \EHAwst \vdash \mForall{\tilde x,\tilde y,\tilde z} \forall{x,y \leq_0 \tilde x,\tilde y} \\
        (\Exists{z \leq_0 t\tilde x\tilde y} \phi(x,y,z) \vee \Forall{z \leq_0 \tilde z} \neg\phi(x,y,z)) &\ \Rightarrow &\text{(logic, arithmetic)} \\
        \EHAwst \vdash \forallst{x,y} \\
        (\Existsst{z \leq_0 txy} \phi(x,y,z) \vee \Forallst z \neg\phi(x,y,z)) &.
      \end{align*}
      If $\EHAwst + \PP \vdash \Phi$, then there is a term~$t$ as above and $t$ induces a computable function that can be used to solve the halting problem, a contradiction.
      \item[$\EHAwst + \PP \nvdash \neg \Phi$.] If $\EHAwst + \PP \vdash \neg\Phi$, then $\EHAwst + \PP$ is inconsistent because $\EHAwst \vdash \neg\neg\Phi$.\qedhere
    \end{description}
  \end{enumerate}
\end{proof}

\section{Transfer}

Nonstandard theories usually include a transfer principle for two reasons:
\begin{enumerate}
  \item transfer connects the standard and the nonstandard universes~\cite{DienerDiener1995};
  \item transfer implies that all uniquely defined objects are standard~\cite[page~1166]{Nelson1977}.
\end{enumerate}
The price to pay for transfer is:
\begin{enumerate}
  \item we can only apply transfer to internal assertions with standard parameters;
  \item transfer may lead to nonconservativity and inconsistency.
\end{enumerate}

The following two principles are inspired by or even taken from Nelson's article~\cite{Nelson1977}.
\begin{definition}
  We define the following principles:
  \begin{enumerate}
    \item the \emph{universal transfer principle}~$\Tforall$ is all instances of $\Forallst{\ul f} (\Forallst{\ul x} \phi \to \Forall{\ul x} \phi)$;
    \item the \emph{existential transfer principle}~$\Texists$ is all instances of $\Forallst{\ul f} (\Exists{\ul x} \phi \to \Existsst{\ul x} \phi)$;
  \end{enumerate}
  where $\ul f$ are all free variables in the internal formula $\phi$.
\end{definition}

Let us see that transfer may lead to nonconservativeness or inconsistency (one solution is to consider instead the rules~$\frac{\Forallst{\ul x} \phi}{\Forall{\ul x} \phi}$ and $\frac{\Exists{\ul x} \phi}{\Existsst{\ul x} \phi}$~\cite[proposition~5.12]{BergBriseidSafarik2012} \cite{Moerdijk1995}). This intuitionistic result (and part of its proof) is an intuitionistic variant of Ferreira and Gaspar's similar classical result (and part of its proof)~\cite[appendix~A]{FerreiraGaspar2015}. We refer to Van den Berg, Briseid and Safarik's article~\cite{BergBriseidSafarik2012} for the notation $\EHAwsst$, $\R$, $\HGMPst$ and $\USzero$ used below.

\begin{theorem}\mbox{}
  \begin{enumerate}
    \item Adding $\Tforall$ or $\Texists$ to $\EHAwsst + \R + \HGMPst$ leads to a nonconservativity over $\HA$.
    \item  Adding $\Tforall$ or $\Texists$ to $\EHAwst$ leads to inconsistency.
  \end{enumerate}
\end{theorem}

\begin{proof}\mbox{}
  \begin{enumerate}
    \item The theory $\EHAwsst + \R + \HGMPst$ proves $\USzero$~\cite[proposition~5.11]{BergBriseidSafarik2012}. The theory $\EHAwsst + \USzero$ plus $\Tforall$ or $\Texists$ is nonconservative over $\HA$~\cite{AvigadHelzner2002} \cite[page~1973]{BergBriseidSafarik2012}.
    \item Let $\neg\st(n^0)$ and $t^1 x \defeq \scriptsize\begin{cases} 0 &\text{if } x <_0 n \\[-0.5ex] 1 &\text{if } x \geq_0 n \end{cases}$. It follows from $\st(\la{x^0}{1})$ and $t \leq^*\! \la{x^0}{1}$ that $\st(t)$. The internal formula $\phi \defequiv tx =_0 0$ is such that $\Forallst x \phi$ and $\neg\Forall x \phi$, contradicting $\Tforall$. The internal formula $\psi \defequiv (tx =_0 1) \wedge \Forall{y < x} (ty =_0 0)$ (notice $\psi \leftrightarrow x =_0 n$) is such that $\Exists x \psi$ and $\neg\Existsst x \psi$, contradicting $\Texists$.\qedhere
  \end{enumerate}
\end{proof}

\section{Recovering standard interpretations}

As mentioned in the introduction, if we restrict our realisability~$\bb$ to the ``purely external fragment'' (where there are only quantifiers of the form $\existsst$ and $\forallst$), then we recover Ferreira and Nunes's bounded modified realisability, and if we restrict our functional interpretation~$\B$ to the ``purely external fragment'', then we recover Ferreira and Oliva's bounded functional interpretation. In this section we sketch a proof of these results, actually on focusing on the second one because it is the more difficult one and the first result is perfectly analogous. Our results and proof are similar to Ferreira and Gaspar's result and proof~\cite[section~4]{FerreiraGaspar2015}. For brevity, we need to assume familiarity with Ferreira and Oliva's theory~$\HAwbd$~\cite[definition~5]{FerreiraOliva2005} and Ferreira and Oliva's bounded functional interpretation~$\BB$~\cite[definition~4]{FerreiraOliva2005}. To state the results, we need to introduce our variant of Ferreira and Gaspar's diamond translation~\cite[page~707]{FerreiraGaspar2015}.

\begin{definition}
  The \emph{diamond translation}~$\diamond$ assigns to each formula $\phi$ of $\HAwbd$ the formula $\phi^\diamond$ of $\EHAwst$ accordingly to the following clauses. For atomic formulas, we define:
  \begin{enumerate}
    \item $(s =_0 t)^\diamond \defequiv s =_0 t$;
    \item $(s \unlhd_\sigma t)^\diamond \defequiv s \leq^\ast_\sigma t$.
  \end{enumerate}
  For the remaining formulas, we define:
  \begin{enumerate}
    \setcounter{enumi}{2}
    \item $(\phi \circ \psi)^\diamond \defequiv \phi^\diamond \circ \psi^\diamond$ for $\circ \in \{\vee,\wedge,\rightarrow\}$;
    \item $(\Finvv{x \unlhd t} \phi)^\diamond \defequiv \Finvv{x \leq^*\! t} \phi^\diamond$ for $\Finv \in \{\Forall,\Exists\}$;
    \item $(\Finvv x \phi)^\diamond \defequiv \Finvvst x \phi^\diamond$ for $\Finv \in \{\Forall,\Exists\}$.
  \end{enumerate}
\end{definition}

Now we state and sketch the proof of our factorisation $\BB\,\diamond = \diamond\,\B$ along the lines of Ferreira and Gaspar's similar factorisation~\cite[proposition~1]{FerreiraGaspar2015}. The factorisation has two formulations: the first one is less clear but (essentially) informs that the terms witnessing $\phi^\BB$ and $(\phi^\diamond)^\B$ are exactly the same, and the second is more clear but less informative.

\begin{theorem}
  For all formulas~$\phi$ of $\HAwbd$, we have:
  \begin{enumerate}
    \item $\EHAwst \vdash \phi_\BB(\ul a;\ul b)^\diamond \leftrightarrow (\phi^\diamond)_\B(\ul a;\ul b)$;
    \item $\EHAwst \vdash (\phi^\BB)^\diamond \leftrightarrow (\phi^\diamond)^\B$.
  \end{enumerate}
\end{theorem}

\begin{proof}\mbox{}
  \begin{enumerate}
    \item The proof is almost a simple induction on the length of $\phi$. For example, in the case of $\exists$ (the most difficult one), we have
    \begin{align*}
      (\Exists x \phi)_\BB(c,\ul a;\ul d)^\diamond &\equiv &\text{(definition of $\BB$)} \\
      (\Exists{x \unlhd c} \mForall{\ul b \unlhd \ul d} \phi_\BB(\ul a;\ul b))^\diamond &\equiv &\text{(definition of $\diamond$)} \\
      \Exists{x \leq^*\! c} \mForall{\ul b \leq^*\! \ul d} \phi_\BB(\ul a;\ul b)^\diamond &\leftrightarrow &\text{(induction hypothesis)} \\
      \Exists{x \leq^*\! c} \mForall{\ul b \leq^*\! \ul d} (\phi^\diamond)_\B(\ul a;\ul b) &\leftrightarrow &\text{(lemma~\ref{lemma interpretation of defined quantifications})} \\
      (\Existsst x \phi^\diamond)_\B(c,\ul a;\ul d) &\equiv &\text{(definition of $\diamond$)} \\
      ((\Exists x \phi)^\diamond)_\B(c,\ul a;\ul d) &.
    \end{align*}

    \item We have
    \begin{align*}
      (\phi^\BB)^\diamond &\equiv &\text{(definition of $\BB$)} \\
      (\mExists{\ul a} \mForall{\ul b} \phi_\BB(\ul a;\ul b))^\diamond &\equiv &\text{(definition of $\diamond$)} \\
      \mExistsst{\ul a} \mForallst{\ul b} \phi_\BB(\ul a;\ul b)^\diamond &\leftrightarrow &\text{(point~1)} \\
      \mExistsst{\ul a} \mForallst{\ul b} (\phi^\diamond)_\B(\ul a;\ul b) &\leftrightarrow &\text{(definition of $\B$)} \\
      (\phi^\diamond)^\B &. &\qedhere
    \end{align*}
  \end{enumerate}
\end{proof}

\bibliography{References}{}
\bibliographystyle{plain}

\end{document}